\documentclass[twoside,leqno]{article}

\usepackage[letterpaper]{geometry}

\usepackage{siamproceedings}

\usepackage[T1]{fontenc}
\usepackage{amssymb}
\usepackage{amsfonts}
\usepackage{graphicx}
\usepackage{epstopdf}
\usepackage{enumitem}
\usepackage{algorithmic}
\ifpdf
  \DeclareGraphicsExtensions{.eps,.pdf,.png,.jpg}
\else
  \DeclareGraphicsExtensions{.eps}
\fi

\newsiamremark{remark}{Remark}
\newsiamremark{hypothesis}{Hypothesis}
\crefname{hypothesis}{Hypothesis}{Hypotheses}
\newsiamthm{claim}{Claim}
\newcommand{\R}{\mathbb R}
\newcommand{\C}{\mathbb C}
\newcommand{\Z}{\mathbb Z}
\newcommand{\N}{\mathbb N}
\newcommand{\T}{\mathbb T}
\newcommand{\e}{\varepsilon}

\newcommand{\cB}{\mathcal{B}}

\newcommand{\cD}{\mathcal{D}}

\newcommand{\cL}{\mathcal{L}}

\newcommand{\cO}{\mathcal{O}}

\newcommand{\di}{{\rm d}}
\newcommand{\pa}{\partial}
\newcommand{\tth}{\mathtt{h}}
\newcommand{\cR}{\mathcal{R}}
\newcommand{\be}{\begin{equation}}
\newcommand{\ee}{\end{equation}}
\newcommand{\uno}{ {\rm Id}}

\newcommand{\grad}{\nabla}

\newcommand{\ii }{{\rm i} }

\newcommand{\om}{\omega}
\newcommand{\vphi}{\varphi}
\newcommand{\im}{{\rm i}}
\newcommand{\tc}{\mathtt{c}}
\newcommand{\Splus}{\mathbb{S}^+} 
\newcommand{\mG}{\mathcal{G}}
\newcommand{\rin}{\mathfrak{r}}
\newcommand{\perd}{\mathtt{d}}
\newcommand{\h}{\mathtt{h}}  
\newcommand{\vect}[2]{\begin{pmatrix}#1 \\#2\end{pmatrix}}
\newcommand{\vet}[2]{\begin{pmatrix}#1 \\ #2 \end{pmatrix}}

\newcommand{\Opbw}[1]{{{{\rm Op}^{{\scriptscriptstyle{\mathrm BW}}}}\big(#1\big)}}
\newcommand{\tm}{\mathtt{m}}
\newcommand{\tG}{{\mathtt G}} 

\newcommand{\ch}{{\mathtt c}_{\mathtt h}}
\newcommand{\ttf}{\mathtt{f}}
\newcommand{\bro}{\bar\rho}
\newcommand{\umu}{\underline{\mu}}
\newcommand{\uphi}{\underline{\varphi}}
\newcommand{\tp}{\mathtt{p}}
\newcommand{\cV}{\mathcal{V}}
\newcommand{\tL}{\mathtt{L}}
\newcommand{\tJ}{\mathtt{J}}

\usepackage{amsopn}

\begin{document}

\newcommand\relatedversion{}

\title{\Large \bf Long time dynamics of space periodic water waves\relatedversion}
    \author{Massimiliano Berti\thanks{SISSA, Via Bonomea 265, Trieste, Italy, \email{berti@sissa.it}, 
    Accepted for the ICM2026 Proceedings to be published by SIAM
 } }
\date{}

\maketitle







\begin{abstract} We review 
recent advances regarding the long-time dynamics of space-periodic water waves, focusing on 
1) bifurcation 
of quasi-periodic solutions, both standing and traveling; 
2) long-time well-posedness 
results; 
3) modulational instability of Stokes waves. These results rely on 
unconventional approaches to KAM and Birkhoff normal form theories for Hamiltonian quasi-linear PDEs 
and  symplectic 
Kato perturbation theory for separated eigenvalues of 
reversible and Hamiltonian operators.  
\end{abstract}



\section*{Introduction.}\label{sec:intro}
The governing equations for water waves originated  in Newton's Principia 
and were formally defined by Euler, Laplace, and Lagrange in the 18th century.
These challenging partial differential equations (PDEs)  
model   a   {\it quasi-linear free boundary} problem and 
have  been  object of extensive research due to their theoretical importance  
and  real-world applications. In recent years, major progress has been made in understanding the long time dynamics of {\it space-periodic}  water waves. 
Previously, global well-posedness results were established for waves with sufficiently 
smooth and localized initial data, whose solutions ``disperse'' and ``scatter'', 
see  \cite{Del}. 
 In contrast,  for periodic initial data, solutions often exhibit oscillatory and recurrent behavior.
 The only known long-term  solutions were {\it Stokes waves}  
--spatially periodic traveling solutions with a speed that depends on the amplitude--
named after Stokes's pioneering  results. 
 Interestingly,  they were the first global in time solutions ever discovered for dispersive PDEs.
Indeed the  water waves equations can be formulated  as a  {\it Hamiltonian} 
quasi-linear dispersive
PDE of the form 
\be \label{diagNL0}
\pa_t u + \ii \Omega( D ) u + N(u, \bar u) = 0 \, ,  \quad D := - \ii \pa_x \, , \quad x \in \T := \R/  2 \pi \Z \, , 
\ee
where $ u(t,x) $ is a complex scalar function,  $ \Omega( D )  $ is a linear Fourier multiplier operator (see \cref{diaglin}) and 
$ N(u, \bar u) $ is a  quadratic 
nonlinearity  
with the same {\it order} of derivatives as the linear term  $ \Omega(D)  $. 

In this Proceeding we provide an overview  of results on space periodic water waves in three key areas:
\\[1mm]
{\bf 1)}  
 {\bf  Kolmogorov-Arnold-Moser (KAM) theory}  \cite{BM,BBHM,BFM1,FG,BFM2}:  
{\it local bifurcation of  ``Cantor-like''  families of 
 time-quasi-periodic water waves,  either standing  or traveling, 
see  \cref{sec:Stokes}.}  

These waves are particularly notable as they provide global-in-time solutions 
with periodic boundary conditions,  for which the initial value problem global 
 well-posedness is still an open question.
\Cref{thm:main00} proves  the existence of  small amplitude quasi-periodic standing waves 
(see 
\cref{fig:StandingW}, left)  
and   \cref{thm:main0} of 
quasi-periodic traveling Stokes waves, which are 
the nonlinear superposition
of  simpler periodic Stokes waves, moving with different non-resonant speeds
(see  \cref{fig:StandingW}, right).
These results are based on a novel approach to KAM theory for $ 1d$ quasi-linear  PDEs. 
\\[1mm]  {\bf 2)}  
{\bf Birkhoff normal  form  (BNF) theory} \cite{BD,BFP,BFF,BMM1}: 
 {\it for any sufficiently smooth and small initial datum of  size $ \varepsilon $,  
the solution  exists and remains of size $ \cO (\e) $, 
for at least times of magnitude $  \varepsilon^{-N} $, see \cref{sec:BNF}.} 

The above KAM quasi-periodic solutions 
exist only for  special Cantor-like sets of space periodic initial data. In contrast, 
these long time existence results  
provide stability information for the solutions starting with  {\it any} initial datum in a sufficiently small 
neighborhood of the origin within a Sobolev space of regular  enough functions.  
In   \cref{teo1}  we describe  the most complete almost global in time stability result known so far
 for periodic water waves. 
These results  rely on a unconventional 
paradifferential Hamiltonian BNF theory for $ 1d$ quasi-linear  PDEs. 
\\[1mm]  {\bf 3)}
 {\bf Modulational instability} \cite{BMV1,BMV3,BMV_ed,BMV4,BCMV}: 
{\it Small amplitude periodic Stokes  
waves are unstable when subjected to wave disturbances with different spatial 
periods, see   \cref{sec:BF}}.

The KAM quasi-periodic solutions 
 are linearly stable under co-periodic space  perturbations, but they may 
become unstable under wave disturbances of different periods.
For periodic Stokes waves  
this phenomenon --known as {\it modulational} or  {\it Benjamin-Feir} instability--
 is now well supported in 
 physical literature, numerical simulations and formal computations, and it may be a key mechanism 
 behind chaotic fluid dynamics.
 Until recently,  there were few rigorous mathematical results. 
\Cref{TeoremoneFinale} and \cref{Thm:1.3}  prove  
new  
existence results of unstable Bloch-Floquet waves 
near  periodic Stokes waves, specifically with long wavelenghts and other 
 periods. 

\smallskip

All these results are particularly challenging because 
the water wave equations are an 
 infinite-dimensional  Hamiltonian dynamical system as in \cref{diagNL0} with a 
{\it non-local}  and  {\it quasi-linear} vector field. 
Even, in the pure gravity case, 
this  is a {\it singular} perturbation problem as  the highest-order 
term of  $ N(u, \bar u) $ in \cref{diagNL0} is a 
nonlinear transport first order operator,  while  $ \Omega(D) = \sqrt{g |D|} $  has 
order $ 1/ 2 $ (cfr. \cref{sec:1.3}). 
\smallskip

KAM and BNF theories were first developed 
for semilinear PDEs,   namely for equations like
\cref{diagNL0} where the nonlinearity $ N(u,\bar u) $  contains no derivatives of $ u $. 
More specifically, 
KAM theory started to be extended around the `80-`90 
 to 
 $1d$  semilinear wave  (NLW)  and Schr\"odinger equations (NLS) by Kuksin,  
Wayne, P\"oeschel, 
Craig, 
Bourgain,  
see e.g. \cite{Kuk}, later extended  in higher dimension 
 in \cite{Bou1,EK}, see also \cite{GYX,BB13,BB20,Pr1} and references therein. 
 We remind that a time quasi-periodic function 
with $ \nu $ frequencies,  with values in a Banach space $ E $,  has the form 
$$
 u (t) = U (\vec  \om t )  
  \qquad {\rm where} \qquad   
  U : \T^\nu \to E \, , \  
  \vphi \mapsto U( \vphi )  \, , \qquad \T^\nu := (\R/2 \pi \Z)^\nu \, , 
$$
is $ 2 \pi $-periodic in the angular variables 
$ \vphi := (\vphi_1, \ldots, \vphi_\nu)  $ and  the frequency  vector $ \vec \om \in \R^\nu $   is 
non-resonant, namely 
$ \vec \om \cdot \ell \neq 0 $ for any $  \ell \in \Z^\nu \setminus \{0\} $. 
In such a case the linear flow $ \{ \vec \om t \}_{t \in \R } $ is {\it dense}  on $ \T^\nu $ and
the torus-manifold 
$  U ( \T^\nu )  $ is invariant under the flow,   which could be ill-posed elsewhere. 
The main challenge  in proving their existence 
is the presence 
 of {\it small divisors}, namely 
terms such as $ \vec \om \cdot \ell $ and other  combinations of 
the normal mode frequencies $ \Omega (j) $, which  accumulate to zero
and appear as denominators in the series expansions.

Regarding Birkhoff normal form, the first results were obtained 
for  semilinear NLS and NLW 
by Bambusi, Gr\'ebert, Bourgain,  Delort,  Szeftel, see \cite{BG03,BDGS} 
and references therein.  
Normal form is a  strategy to  extend the lifespan of solutions of PDEs. 
First applied for semilinear wave equations in \cite{S84}, the core idea is to construct 
changes of variables that iteratively reduce the {\it size} of the nonlinearity.
However   ``resonant''  terms,  corresponding to a combination of normal mode frequencies 
$ \Omega (j_1) \pm  \ldots \pm \Omega (j_N)   = 0 $, 
cannot be eliminated and remain in the normal form.  
In order to guarantee that the resonant normal form does not cause the 
explosion of solutions, 
the Hamiltonian or reversible structure of the vector field plays a crucial role. 
In addition, too small   $N$-{\it wave resonances}  
$ \Omega (j_1) \pm  \ldots \pm \Omega (j_N)  $ 
can cause the normal form transformations to be unbounded, and must  be avoided, as achieved
in 
\cite{BG03,BDGS}.

After these pioneering results a lot of effort 
 has been  devoted  
 to understanding the effect of derivatives in the  nonlinearity 
 in KAM and BNF theory, see e.g. \cite{Kuk,KP02,LY,BBP1,BBP2} for semilinear PDEs 
 with derivatives.  
 The primary goal was 
 to determine if such theories are applicable to quasi-linear PDEs of fundamental physical importance, 
like  the water wave and Euler equations.
All the previous  approaches are  not applicable
since  they merely produce  formal, unbounded changes of variables.

\smallskip

For quasi-linear PDEs, 
the KAM and BNF strategies required a fundamental shift in perspective.
The key criterion for what is considered ``small'' in the nonlinearity $ N(u, \bar u) $  
 is not the \underline{size} of the nonlinearity itself (quadratic, cubic, etc,  in $u $), 
 but rather the \underline{order} of the derivatives of $ u $ it contains. 
This change in paradigm 
was accomplished by systematically introducing
 {\it pseudo differential}  calculus and  {\it Fourier integral operators} ideas 
in KAM theory for PDEs and their 
nonlinear paradifferential analogue in BNF theory, as we shall explain in 
\cref{sec:Stokes,sec:BNF}. 
The new approaches to KAM and BNF theories 
may be presented 
as follows. 
\\[1mm]
 {\bf 1)}  {\bf Pseudo-differential KAM normal form}. 
Pseudo differential calculus and  Fourier integral operators  
are used to {\it reduce}  the  unbounded linearized equations at an  approximate quasi-periodic  solution
along the Nash-Moser algorithm,  to {\it constant coefficients} Fourier multiplier operators,   
up to {\it smoothing} remainders, 
see \cref{sec:1.3}, specifically comment (a) below \cref{LINOP}. 
Ideas of this kind were pioneered in \cite{PlTo,IPT} for time periodic solutions. 
A  KAM argument then completes 
the full diagonalization of the equations, 
providing a  
very precise high-frequency asymptotic expansion 
of  the perturbed eigenvalues, see comment (b) below \cref{LINOP}.  
\\[1mm]
 {\bf 2)}  {\bf Para-differential BNF normal form}.  Paradifferential calculus 
--a nonlinear version of pseudodifferential calculus-- 
 allows for an  analogous reduction of nonlinear PDEs to {\it constant coefficient} symbols, modulo 
 {\it smoothing remainders}, see \cref{sec:ideasBNF}.
 This procedure is  not symplectic, but  we may systematically 
 restore the Hamiltonian structure of the equation, at any  
  degree of homogeneity, by a symplectic 
  Darboux perturbative argument. 
 At this point we  implement BNF transformations to decrease the size of 
 the symbols and of the smoothing remainders, see comments (a)-(c) below \cref{equaz1}.
 Another Hamiltonian paradifferential BNF approach 
specific to quasi-linear Klein-Gordon  equations  was  developed in   \cite{Del1,Del3}.

\smallskip

The modulational instability of Stokes waves in deep water was an unexpected discovery made by Benjamin and Feir in 1967. 
They were trying to create periodic waves in a tank, but, no matter how carefully the
experiment was set up, the waves always disintegrated traveling down the canal. 
Recognizing the 
physical significance of this issue,  
Benjamin and Feir  \cite{BF} proposed a heuristic justification 
of the instability, based on  energy transfer between waves
of nearby frequencies (side-band instability). 
Related  explanations  were independently advanced 
by  Whitham, Zakharov and Lighthill  for also other physical phenomena.
The first rigorous mathematical results  for Stokes waves were established in 
 \cite{BrM} 
 via spatial dynamics and center manifold theory 
 and recently in \cite{NS} via a Lyapunov-Schmidt decomposition. 
 A more complete mathematical description of such instabilities has been achieved, 
and is rapidly progressing  through a novel spectral approach, initiated in  \cite{BMV1}. Its 
  core idea is a  
\\[1mm]
 {\bf 3)}  {\bf Symplectic and reversible Kato perturbation theory}
 to analyze   the splitting of multiple 
separated eigenvalues of the linearized water waves operators at the 
Stokes wave, 
as the 
spatial period varies (Floquet parameter),  
combined with block-diagonalization algorithms inspired by KAM theory, see \cref{sec:ideasBF}.

\smallskip

In all  these results the {\it Hamiltonian}/{\it reversible} structure of the 
water waves equations --previously unexploited in rigorous analysis--
plays a key role.
For the global well-posedness results  
 \cite{Wu2,GMS,IP,DIPP,AlDe1}  with sufficiently 
smooth, small, and localized in space  initial data in $ \R^d $,    these properties 
are not relevant,  
because the dispersive properties of the flow are dominant. 
If $ d= 2 $, 
the water waves solutions asymptotically scatter  
to a  free linear wave. If $ d= 1 $, the dispersion is weaker, leading to modified scattering
of the  solutions  
to some linear wave with a small nonlinear adjustment.
 In contrast, for space-periodic  water waves, there is no dispersion and 
 the nonlinear long-term behavior 
 is heavily influenced by $N $-wave resonant interactions and the Hamiltonian and reversible nature of the vector field.  

\smallskip

The paper is structured  as follows. 
We begin with a precise introduction to the water waves equations and their properties. Then in  \cref{sec:Stokes,sec:BNF,sec:BF}
we present the results in detail 
along with  related literature and ideas of proofs. 
Each section is concluded with a list of challenging open problems,
that warrant particular attention. 
\subsection{The water waves equations.} 
\begin{figure}[h]
    \begin{minipage}[c]{0.5\textwidth}
        \centering
        \includegraphics[width=8cm]{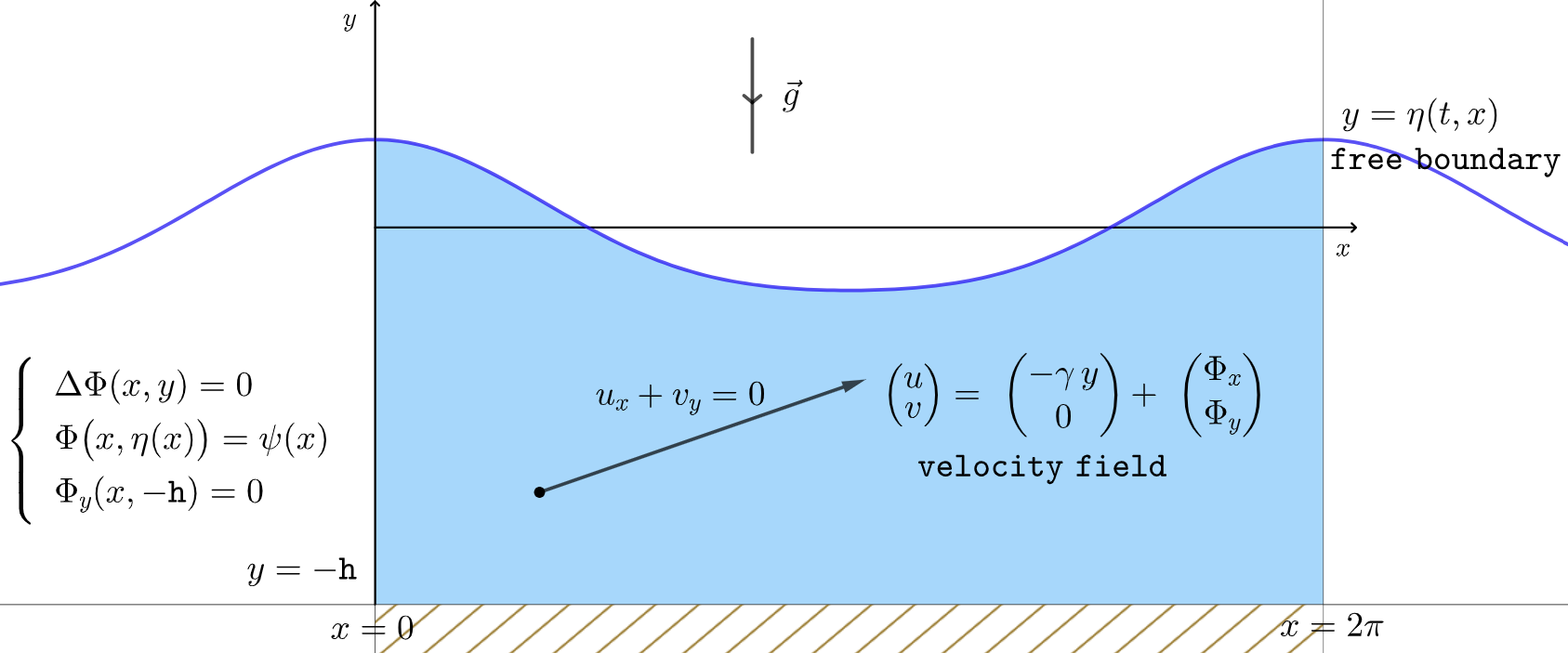}
        \label{fig:waterwaves}
    \end{minipage}\hfill
    \begin{minipage}[c]{0.5\textwidth}
        \centering
        \caption{The gravity water waves problem. 
The velocity field inside the time dependent domain $ \cD_\eta $ is divergence free and with constant vorticity $ \gamma $. Capillary forces may act at the free interface $ y = \eta (t,x) $.}
    \end{minipage}
\end{figure}

We consider the time evolution of a 2D
inviscid and incompressible 
fluid with constant vorticity, which occupies  
the time dependent region 
\be\label{eq:domain}
\cD_{\eta} := \big\{ (x,y)\in \T\times \R \ : \ - \tth \leq y \leq \eta(t,x) \big\} \, ,
\ee
delimited by a bottom  with depth $ \tth > 0 $, possibly infinite, 
and space periodic boundary conditions. The fluid's movement is driven 
by gravity and possibly capillary forces at the free surface. 
The dynamical unknowns are the  free boundary  $ y = \eta (t, x)$ of $\cD_{\eta} $ 
and the divergence-free  velocity field
$ \big( \begin{smallmatrix} u(t,x,y) \\ v(t,x,y) \end{smallmatrix} \big) $, 
whose vorticity $ v_x - u_y = \gamma $ is constant in $ \cD_{\eta} $, 
a property  preserved throughout the time evolution.
In view of the Hodge decomposition,  the velocity field  
is the sum of the Couette flow $ \big( \begin{smallmatrix} - \gamma y \\ 0 \end{smallmatrix} \big) $, which carries all the  vorticity $ \gamma $, and an irrotational field 
$ (\Phi_x,\Phi_y)  $, where $ \Phi (t, x,y) $  is called the generalized velocity potential, see  \cref{fig:waterwaves}.

We study the water waves problem in the Zakharov-Craig-Sulem Hamiltonian formulation
 \cite{Zak1,CrSu} extended  for constant vorticity fluids in  \cite{CIP,Wh}.  
Denoting by 
$\psi(t,x)$  the evaluation of the generalized velocity potential at the free interface
$ \psi (t,x) := \Phi (t,x, \eta(t,x)) $, one recovers
$ \Phi $ 
as the unique harmonic solution of 
$$  
 \Delta \Phi (t,x,y) = 0 \  \text{in}  \ {\mathcal D}_\eta \, , \quad  
  \Phi(t,x,\eta(t,x)) = \psi(t,x) \, , \quad 
 \Phi_y(t,x, - \tth ) = 0  \, , 
$$
or  $\Phi_y ( t, x, y) \to  0  $ as $ y \to  - \infty $ in deep water. 
The latter Neumann boundary condition 
 means the impermeability property of the bottom  $ y = - \tth $ of $ \cD_\eta $. 

Inside $ \cD_{\eta}$ the velocity field of the fluid evolves according to
the Euler equations and the 
dynamics of the free surface is  
determined by 
two boundary conditions. 
The first  is 
that the fluid particles on the time dependent surface remain on it along the evolution
(kinematic boundary condition).
The second condition is that   the pressure of the fluid 
plus the gravity and capillary forces at the free surface  is equal to the constant 
atmospheric pressure (dynamic boundary condition). 
It turns out that $ (\eta(t,x), \psi (t,x)) $ evolve according to the 
quasi-linear  equations 
\begin{equation} \label{eq:etapsi}
\left\{\begin{aligned}
 &   \partial_t \eta = G(\eta)\psi+ \gamma \eta \eta_x \\
&\partial_t\psi =  -g\eta  -\frac{1}{2} \psi_x^2 
+  \frac{1}{2}\frac{(\eta_x  \psi_x + G(\eta)\psi)^2}{1+\eta_x^2}+\kappa 
\partial_x\Big[ \frac{\eta_x}{(1+\eta_x^{2})^{\frac{1}{2}}}\Big]+ \gamma \eta \psi_x+ \gamma \partial_x^{-1} G(\eta) \psi 
\end{aligned}\right.
\end{equation}
where $ g > 0 $ is the gravity constant, $ \kappa \geq 0 $ is the surface tension coefficient,
$ \partial_x\big[ \frac{\eta_x}{(1+\eta_x^{2})^{1/2}}\big] $
is the curvature of the  surface, $ \gamma $ is the constant vorticity 
and  
$$ 
G(\eta)\psi  := (- \Phi_x \eta_x + \Phi_y)\vert_{y = \eta(t,x)}  
$$ 
is the  Dirichlet-Neumann operator. The name refers to the fact that  $G(\eta) $ 
transforms the Dirichlet datum  $ \psi$ 
of the  velocity potential $ \Phi $ into its Neumann derivative $G(\eta)  \psi $ at the free surface.  
 Such operator has been  extensively investigated, see e.g. \cite{LannesLivre}. 
It 
 is linear in $ \psi $, self-adjoint and non-negative,  
it vanishes only  on the constants, i.e.  $ G(\eta)[1] = 0 $.
The operator $ G(\eta) $ depends  analytically on the boundary $ \eta (x) $ in 
several  topologies. 
Furthermore $G(\eta) $  is a 
pseudo-differential operator,  if $ \eta (x) $ is $ C^\infty $.
For  2D fluids, it has the asymptotic expansion 
$ G(\eta) = D \tanh (\tth D) + \cR (\eta) $ where $  D = - \ii \pa_x  $ 
and 
$ \cR (\eta) $   is an infinitely many times regularizing operator, that vanishes for $ \eta = 0 $. 
\\[1mm]
{\sc Phase spaces.}
The space average  $ \langle \eta \rangle := 
\frac{1}{2\pi } \int_\T \eta(x) \, \di x $ is a first integral of \cref{eq:etapsi}.
Additionally 
the vector field on the right hand side of \cref{eq:etapsi} depends only on $ \eta $ and 
$ 
 \psi - \langle \psi \rangle   $. As a consequence, 
the variables $ (\eta, \psi) $ of system \cref{eq:etapsi} can be set to belong to 
$ H^s_0(\T) \times \dot H^s (\T) $.  
Here  $H^s_0(\T)$ 
 denotes the Sobolev space of zero average functions 
and  $\dot H^s(\T)$    the homogeneous space where functions 
that differ by a constant are considered equivalent.
\\[1mm]
{\sc Hamiltonian formulation.} As observed in \cite{Zak1,CrSu,CIP,Wh} 
the equations \cref{eq:etapsi} are the Hamiltonian system 
\be
\begin{aligned}
& 
\qquad \pa_t 
\begin{pmatrix} 
\eta \\ 
\psi
\end{pmatrix} 
= J_\gamma \vect{\nabla_{\eta} H(\eta,\psi)}{\nabla_{\psi} H(\eta,\psi)} \qquad
\text{where} 
\qquad 
J_\gamma := \begin{pmatrix} 0 & \uno \\ -\uno & \gamma \pa_x^{-1} \end{pmatrix} \qquad \text{and} \\
& \label{H.gamma}
H (\eta, \psi)  := 
\underbrace{\frac12 \int_\T 
\big(\psi \, G(\eta ) \psi+ g \eta^2\big) \, \di x}_{kinetic  + potential \ energy}  
+ \underbrace{\kappa  
\int_{\T}  \sqrt{1+ \eta_x^2}  \, \di x}_{capillary \ energy}+  
\underbrace{\frac{\gamma}{2} \int_{\T}  \big(-\psi_x \eta^2+ \frac{\gamma}{3} \eta^3\big)  \, \di x}_{vorticity \ energy} \, .
\end{aligned}
\ee
The $L^2$-gradients $(\grad_\eta H, \grad_\psi H) $ in
\cref{H.gamma}  belong to (a dense subspace of)  
$ \dot L^2(\T)\times L^2_0(\T)$. The Poisson tensor $ J_\gamma $ is not the canonical one for $ \gamma \neq 0 $, but Wahl\'en noted in  \cite{Wh} that
the Hamiltonian system \cref{H.gamma} assumes the 
standard Darboux  form in the 
canonical variables 
$ (\eta, \zeta) := (\eta, \psi -  \frac{\gamma}{2} \pa_x^{-1}  \eta ) $. 

\begin{figure}\label{fig:StandingW}
 \centering
\includegraphics[width=13cm]{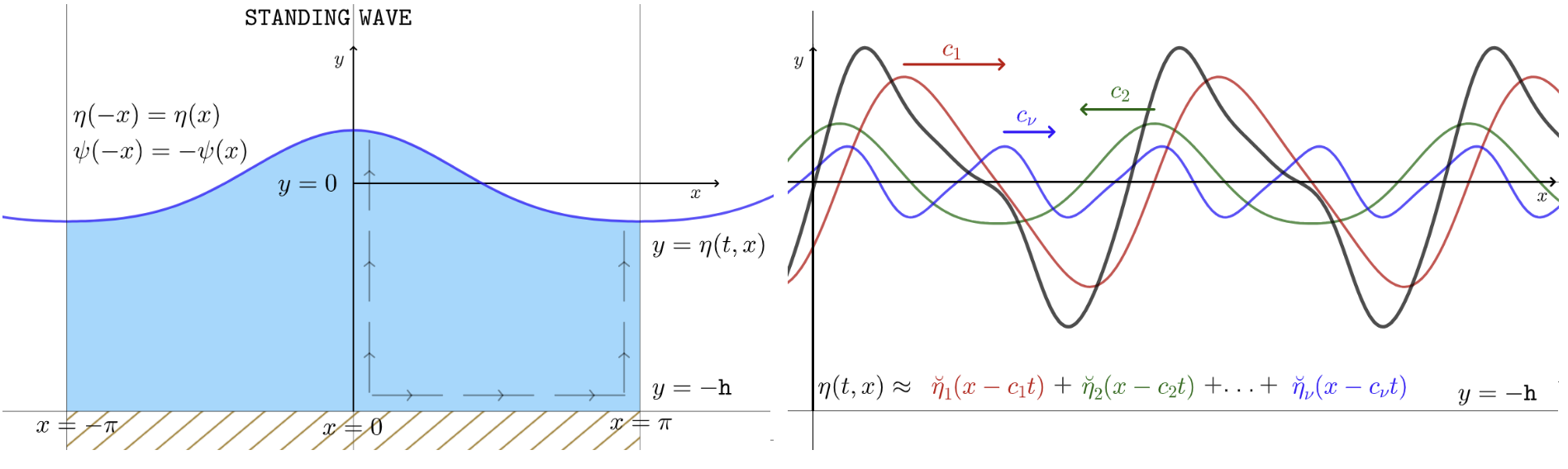}
\caption{ On the left, standing waves for irrotational fluids ($\gamma =0 $); these solutions are even in $ x $ at any time. 
On the right, quasi-periodic traveling Stokes waves with $ \nu $ frequencies, see \cref{def:QPT}.}
\end{figure}

\noindent 
{\sc Reversible structure.} In addition to being  Hamiltonian,
the water wave system
  \cref{eq:etapsi} 
is time-reversible with respect to the involution
$$
\rho\vet{\eta(x)}{\psi(x)} := \vet{\eta(-x)}{-\psi(-x)} \, , \quad 
\text{i.e.} \  H \circ \rho = H \, .
$$
{\sc Space invariance.}  Furthermore, 
because the bottom of the fluid domain $ \cD_\eta $ in \cref{eq:domain} 
is flat, system \cref{eq:etapsi}  is 
invariant under space translations: 
 the Hamiltonian $ H  $   is preserved by the symplectic maps 
  $$ 
  \tau_\varsigma : u(x) \mapsto 
  u(x + \varsigma ) \, , 	 \quad u(x) = \vet{\eta(x)}{\psi(x)} \, , \quad 
  \text{namely} \quad  
  H \circ \tau_\varsigma = H \, , \  
  \quad \forall \varsigma \in \R \,  . 
 $$
   This symmetry
   is essential for  the existence of 
   traveling waves.
\\[1mm] 
{\sc Standing waves.} For  irrotational water waves, i.e. for  $ \gamma = 0 $, 
the equations \cref{eq:etapsi} have  another important symmetry:  
the subspace of functions  $ (\eta(x), \psi(x)) $ 
which are even in $ x $  --called {\it standing waves}-- 
remains invariant under time evolution. 
This is not the case when  $ \gamma \neq 0 $. 
This symmetry implies that $ \Phi_x (0,y) = \Phi_x (\pi,y) = 0 $, namely 
the velocity field $   (u,v) = (\Phi_x, \Phi_y ) $ is purely vertical at $ x = 0  $ and 
$ x = \pi $.
This condition is equivalent to the physical problem of water waves confined 
in a tank with impermeable vertical walls, see \cref{fig:StandingW} left. 

\smallskip

We describe the dynamics of system  \cref{eq:etapsi} near the equilibrium solution $ (\eta,\psi) = (0,0) $, 
where the ocean is flat and still.  

\subsection{Linear water waves.} 
The linearized equations  \cref{eq:etapsi} at the equilibrium $ (\eta,\psi)=(0,0)$ are
\be\label{eq:lin}
\pa_t \vect{  \eta}{  \psi}  = 	
\begin{pmatrix} 0  & \!\!\! G(0) \\ -(g+ \kappa D^2  ) & 
 \gamma  \pa_x^{-1} G(0)\end{pmatrix} \!\!\vect{ \eta}{  \psi}  \, , 
\qquad G(0) = D \tanh (\tth D)\, , \qquad D = \frac{1}{\ii} \pa_x \, , 
\ee
where  
$G(0)  $ 
is the Fourier multiplier 
with  symbol
$  \tG(\xi):=   \xi\tanh(\tth \xi) $ if $ 0<\tth<+\infty $ and  
$   \tG(\xi):=  |\xi| $ if $  \tth=+\infty $.

System \cref{eq:lin} may be diagonalized  introducing the complex variable 
\be\label{defMM-1}
u = 
\frac{1}{\sqrt{2}} (M^{-1}(D)  \eta    +   \ii M(D)  \big( 
\psi -  \tfrac{\gamma}{2} \pa_x^{-1}  \eta ) \big)
\quad
\text{where} \quad  
M(D):= \Big( \frac{G(0)}{ g+ \kappa D^2 + \frac{\gamma^2}{4} G(0)D^{-2} } \Big)^{1/4} \, . 
\ee
A direct computation
 shows that 
\be \label{diaglin}
\pa_t u = -\ii \Omega(D) u \quad \text{where} \quad  
\Omega(D):= \sqrt{ G(0)
 \Big( g+ \kappa D^2 + \frac{\gamma^2}{4} \frac{G(0)}{D^{2}} \Big)} + \ii \frac{ \gamma}{2} G(0) \pa_x^{-1}   
\, .
\ee
If $ \kappa > 0  $ 
the dispersion relation $ \Omega (D) \sim 
\sqrt{\kappa}|D|^{\frac32}  $ grows super-linearly whereas in the pure gravity case
 $ \Omega (D) \sim \sqrt{g}|D|^{\frac12} $.
In the Fourier expansion $u(x) = 
\sum_{j \in \Z\setminus \{0\}} u_j \,  e^{\im j x} $, equation \cref{diaglin} amounts 
to  infinitely many harmonic oscillators
$ \pa_t u_j = - \im \Omega (j) u_j $, $  j \in \Z \setminus \{0\} $, 
with {\it normal mode} 
frequencies of oscillations 
\be\label{omegonejin}
\Omega (j) := 
\sqrt{ \tG(j) 
\Big(g+ \kappa j^2 + \frac{\gamma^2}{4} \frac{\tG(j)}{j^{2}}\Big)}  +  \frac{ \gamma}{2} \frac{\tG(j)}{j}  \, ,
\quad \tG(j) = j \tanh (\tth j) \, , 
\ee
The equilibrium  $ (\eta,\psi) = (0,0) $ is thus elliptic. 
All the solutions of \cref{diaglin} are the 
superposition of {\it plane waves} 
\be\label{eq:solline}
u(t,x) = {\mathop \sum}_{j \in \Z \setminus \{0\} } u_j(0) e^{\ii \Omega (j) 
 t } e^{\ii j x }
\ee
each having wave speed  $ \Omega (j) $ and wave vector $ j $. 
The linear solutions \cref{eq:solline} 
are periodic, quasi-periodic, or almost periodic in time according to 
the irrationality properties of the frequencies 
$ \{ \Omega (j)  \}_{j \in \Z \setminus \{0\} } $, that 
depend on the   {\it physical parameters} $ g, \kappa, \gamma, \tth $, 
and the number of non zero initial amplitudes $ u_j (0)  $.

\subsection{Major questions.}
In the following \cref{sec:Stokes,sec:BNF,sec:BF} we will address these key questions:
\begin{enumerate}
\item[{\bf Q1)}] Do quasi-periodic solutions,  either standing or traveling, 
of  the nonlinear water wave equations exist? 
\item[{\bf Q2)}] What is the maximum time interval of existence of small-amplitude space periodic 
water waves solutions?
\item[{\bf Q3)}] Are the periodic  Stokes water waves  stable or unstable under wave disturbances of different spatial periods? 
\end{enumerate}

The primary difficulties in answering {\bf Q1)}-{\bf Q3)} stem from the  quasi-linearity
of the water waves vector field.

\section{Quasi-periodic water waves.}\label{sec:Stokes}

In 1847 
Stokes  discovered  
the existence of pure gravity spatially periodic traveling water waves 
of the form 
\begin{equation}\label{eq:Stokes1}
\eta(t,x)=\breve \eta(x-ct) \, , 
\quad \psi(t,x)= \breve \psi(x-ct) \,  ,
\end{equation}
 for some speed  $c$
 and  $2\pi$-periodic profiles  $\breve \eta (x), \breve \psi (x) $. Such waves are $ 2 \pi / c $ 
 time periodic  and persist indefinitely despite dispersive effects. 
 Their first mathematically rigorous existence 
 proof  was given by Levi-Civita and Nekrasov: 
there is  a unique family  of real analytic 
functions $( \breve \eta_\e(x), \breve \psi_\e(x), c_\e)$, parameterized by the amplitude $|\e| \leq \e_0$,  such that
 $ \breve \eta_\e(x-c_\e t), \breve \psi_\e(x - c_\e t) $ are solutions of  
the pure gravity water waves equations \cref{eq:etapsi} with 
$ \kappa = 0 $ and $ \gamma = 0 $. 
The profiles 
 $ \breve \eta_\e (x), \breve \psi_\e (x) $ are $2\pi$-periodic, 
 and,  for gravity  $ g = 1 $,  
 they have  the Taylor  expansion  
\begin{equation}\label{eq:Stokes2}
 \begin{aligned}
&  \breve \eta_\e (x) =  \e \cos (x) 
+ 
  \cO(\e^2) \, , \ 
  \breve \psi_\e (x) = \e 
  \tc_\tth^{-1} \sin (x) 
  +\cO(\e^2) \, , 
  \ c_\e = 
  \tc_\tth + \cO(\e^2) \, ,
 \end{aligned}
\end{equation}
where 
$ \tc_\tth :=  \sqrt{\tanh(\tth)} 
     $ and 
     $ \tc_\tth = 1 $ 
     if $ \tth = + \infty $. 
The  literature concerning Stokes waves 
is currently enormous and we refer to the recent survey \cite{ebb} for an updated presentation.

\subsection{Quasi-periodic standing waves.}
Building on previous results for 
quasi-linear KdV and NLS equations 
 \cite{BBM14,BBM16,FP}, the first  quasi-periodic water waves were constructed
  in the gravity-capillary case in \cite{BM} and later in the pure gravity case in \cite{BBHM}, 
 demonstrating the existence of standing  irrotational 
 waves. Earlier constructions of time-periodic standing waves were proved in
  \cite{PlTo, IPT}  for pure gravity waves and  in  \cite{AB} in the  gravity-capillary case.
  We denote $ \R_+ $ the real positive numbers.

\begin{theorem} \label{thm:main00}  {\bf (Quasi-periodic standing water waves \cite{BBHM})}
Consider finitely many ``tangential'' sites $ \Splus \subset \N  $
and denote $ \nu := | \Splus| $ its cardinality. 
Fix a subset 
	$ [\tth_1, \tth_2 ] \subset (0,+\infty) $.   
Then there exists $ \bar s > 0 $,  
$ \e_0 \in (0,1) $ such that 
for every $ |\xi |   \leq \e_0^2  $, $  \xi := (\xi_j)_{j \in \Splus }  \in \R_+^\nu $, 
the following holds: 
\begin{enumerate}
\item there exists 
a Borel set  $ \mG_\xi \subset [{\mathtt h}_1, {\mathtt h}_2] $ 
with density one at $ \xi = 0 $, i.e. 
$ \lim_{\xi \to 0} | \mG_\xi  |  = {\mathtt h}_2- {\mathtt h}_1 $;
\item 
for any  $ {\mathtt h} \in \mG_\xi  $,  the pure gravity irrotational 
water waves equations \cref{eq:etapsi}
(with zero surface tension $ \kappa = 0 $ and  vorticity $ \gamma = 0 $)
have a time quasi-periodic standing wave solution 
$ 
 (\eta ( \vec \om t, x), \psi (  \vec\om t, x) ) $ 
with Sobolev regularity $ (\eta, \psi)  (\varphi,x) 
  \in H^{\bar s} ( \T^\nu \times \T, \R^2) $,    of the form
\be\label{QP:soluz0}
\begin{aligned}
& \eta(\vec \omega t,x) = {\mathop \sum}_{ j \in \Splus}  \sqrt{\xi_j} \cos ({ \om}_j t)  \cos(jx) + r_1 (  \vec  \om t,x ) \, , \\
& \psi(\vec\omega t ,x) = - {\mathop \sum}_{ j \in \Splus} \sqrt{\xi_j}  \Omega(j)^{-1}  \sin ({  \omega}_j t)  \cos(jx)+ r_2 ( \vec  \om t,x ) \, , 
\end{aligned}
\ee
where   $ \Omega(j) = |j|^{\frac12} \tanh^{\frac12}(\mathtt h |j|) $ are 
  the unperturbed linear normal mode frequencies, 
\begin{enumerate}
\item the remainders
$ r_i ( \vphi ,x ) $  satisfy 
  $ \lim_{\xi \to 0} \frac{ \| r_i \|_{H^{\bar s}}}{\sqrt{|\xi|}} = 0 $  for $ i = 1, 2 $;
\item the frequency vector $ \vec  \omega := ( \omega_j)_{j \in \Splus} $ depends on  $\xi $,
is Diophantine  and  satisfies
$ \lim_{\xi \to 0}\vec { \omega} = 
(\Omega(j))_{j \in {\mathbb S}^+} $. 
\end{enumerate}
\end{enumerate}
In addition the quasi-periodic solutions \cref{QP:soluz0} are  linearly stable.  
\end{theorem}

Let us make some comments.  
\\[1mm]
1.  {\sc Linear stability.} We prove that the linearized water waves equations \cref{eq:etapsi}  at the 
quasi-periodic solutions \cref{QP:soluz0},
are conjugated, in the  directions normal  to the invariant torus,  to infinitely many 
scalar autonomous ODEs 
\be\label{eq:linKAM}
\partial_t {\mathtt w}_{j} = - \ii \mu_j {\mathtt w}_{j} \, , 
\quad \forall j \in  \Z \setminus \mathbb S_0 
\quad \text{where}  \quad {\mathbb S}_0 := \Splus \cup (- \Splus) \cup \{ 0 \} 
\ee
and 
$$
 \mu_j 
:= \mathtt m_{\frac12} |j|^{\frac12} \tanh^{\frac12}(\mathtt h |j|) 
+ \rin_j \in \R \, , 
\quad \rin_j = \rin_{-j} \,, \quad 
\mathtt m_{\frac12} = 1 + \cO( | \xi |^{c}) \, , \quad 
\sup_{j \in  \Z \setminus \mathbb S_0 } |j|^{\frac12} |\rin_j| = \cO( | \xi |^{c} ) \, , 
$$ 
for some $ c  >  0 $. 
The $ \ii \mu_j $ are the non zero 
 {\it Floquet exponents} of the quasi-periodic solutions. 
 This is not only a dynamically relevant information 
but also a key ingredient in proving the existence of quasi-periodic solutions,
as we comment below.  
Note that for $ \xi = 0 $ the $  \mu_j $  coincide with the unperturbed linear normal mode 
frequencies of  oscillations $ \Omega(j) = |j|^{\frac12} \tanh^{\frac12}(\mathtt h |j|)  $. The multiplicative 
contant $ \mathtt m_{\frac12}   $
takes into account the quasi-linear  effects of the nonlinearity.  
\\[1mm] 2. 
{\sc Melnikov non-resonance  conditions.} 
In  \cref{thm:main00} 
we eliminate a small measure set of depths $ \tth   $ 
to impose non-resonance conditions between normal mode frequencies. We require that
 the frequency vector $ \vec \omega \in \R^\nu $ is  {\it diophantine}, 
namely there exist $ \upsilon >  0 $ and $ \tau > \nu - 1 $ such that   
\be \label{dioph}
| \vec \om \cdot \ell | \geq \upsilon |\ell|^{-\tau} \, , \quad \forall \ell \in \Z^{\nu} \setminus \{0\} \, , 
\ee 
and, additionally,  
  first and second Melnikov non-resonance conditions as 
\begin{align}\label{1nd melnikov perd}
&  | \vec \omega \cdot \ell  + \mu_j 
   | \geq  \upsilon j^{\frac12} \langle \ell  \rangle^{- \tau}, \, 
 \forall \ell   \in \Z^\nu, \, j \in \N \setminus \Splus, \\
& \label{2nd melnikov perd}
 |\vec \omega \cdot \ell  + 
 \mu_j -  \mu_{j'} 
 | \geq 
  \upsilon j^{-\perd} j'^{-\perd} \langle \ell  \rangle^{-\tau}, \,\,
 \forall \ell   \in \Z^\nu ,\,\,j, j' \in \N \setminus \Splus , \, (\ell,j,j') \neq (0,j,j) \, ,
\end{align}
for some $ \perd > 0 $. 
We  can impose only weak non-resonance conditions such as \cref{2nd melnikov perd}, 
which lead to a loss of $\cO(\perd) $ space derivatives, because the linear frequencies
$ \mu_j  \sim \sqrt{j}  $ grow only sublinearly as $j \to + \infty $. 
In  contrast, all previous KAM theorems, see e.g.  \cite{Kuk}, 
required that the eigenvalues of the linear unperturbed 
differential operator grow 
as $ j^d $ with $ d \geq 1 $.
We compensate for this loss of derivatives using a regularization procedure described 
below \cref{LINOP}.

\subsection{Quasi-periodic  Stokes  waves.}
We now present an 
existence result of  quasi-periodic traveling 
pure gravity water waves, i.e. for 
system \cref{eq:etapsi} with zero surface tension $ \kappa = 0 $, and with possibly non-zero 
vorticity $ \gamma $.  We set $ g = 1 $ 
so  the linear frequencies are $  \Omega (j)  = 
	\sqrt{ j\tanh (\tth j ) \big(1+  \frac{\gamma^2}{4} \frac{\tanh (\tth j )}{j}\big)}  +
	 \tfrac{\gamma}{2} \tanh (\tth j) $.

\begin{definition}  {\bf (Quasi-periodic traveling wave)} \label{def:QPT} 
	{\it We call $ (\eta (t,x), \psi(t,x)) $  a 
	time quasi-periodic {\it traveling} wave  
	with non-resonant frequency vector  $ \vec \omega = ( \omega_1, \ldots, \omega_\nu)  \in \R^\nu $, 
	and ``wave vector'' $ \vec \jmath 
	= ( j_1, \ldots, j_\nu)  \in \Z^\nu $,  if there exist 
	functions
	$ ( \breve\eta, \breve\psi ) : \T^\nu \to \R^2 $ such that} 
	$ (\eta ( t, x), \psi ( t, x))  
	= ( \breve\eta( \vec \omega  t - \vec \jmath  x), 
	\breve\psi ( \vec \omega  t - \vec \jmath  x))  $. 
\end{definition}
If $ \nu = 1 $ they have the form \cref{eq:Stokes1}.  
If $ \nu \geq 2 $ 
solutions of this kind  never look stationary in any moving frame.

Fixed   finitely  many  arbitrary  
integers 
$ \bar n_1, \ldots, \bar n_\nu   $ satisfying 
	$ 1 \leq \bar n_1 < \ldots < \bar n_\nu $  and {\it signs}   
	$  \sigma_1 , \ldots , \sigma_\nu \subset 
	 \{ -1, 1 \} $, 
we consider the quasi-periodic traveling wave solutions   of the linearized system 
\cref{eq:lin} given by 
$$
\begin{pmatrix}
	        \breve \eta_L (\vec \Omega  t - \vec \jmath x) 
	        \\ \breve \psi_L (\vec \Omega  t - \vec \jmath x)
		\end{pmatrix} \  
		\text{with frequency vector} \
		\vec \Omega  := ( \Omega (\sigma_a \bar n_a) )_{a=1, \ldots, \nu} \ 
		 \text{and wave vector} \ 
	  \vec \jmath := (\sigma_a \bar n_a )_{a=1, \ldots,\nu} \, , 
$$
where, for $ \vartheta := (\vartheta_a)_{a=1,\ldots,\nu} \in \T^\nu $,  
\be\label{sel.sol}
 \begin{pmatrix}
	        \breve \eta_L (\vartheta) \\ \breve \psi_L (\vartheta)
		\end{pmatrix}  :=  \!\!  \!    \sum_{a \in \{1, \ldots, \nu \colon   \sigma_a = + 1\}}   \!  
		\begin{pmatrix}
			M_{\bar n_a} \sqrt{\xi_{\bar n_a}} \cos ( \vartheta_a)  
			\\ 
			P_{\bar n_a} \sqrt{\xi_{\bar n_a}} \sin ( \vartheta_a) 
		\end{pmatrix}     +   \!\! \!     \sum_{a \in \{1, \ldots, \nu \colon   \sigma_a = - 1\}}
		\! \! \begin{pmatrix}
			M_{- \bar n_a} \sqrt{\xi_{- \bar n_a}} \cos ( \vartheta_a) \\ 
			P_{-\bar n_a} \sqrt{\xi_{- \bar n_a}} \sin ( \vartheta_a) 
		\end{pmatrix} \, .
		\ee
Here $ \xi_{\pm \bar n_a} >  0 $, 
$
		M_j := \big( \tfrac{\tG(j)}{g +
			\frac{\gamma^2}{4}  \tG(j) j^{-2}} \big)^{\frac14} $,  
		$ j \in \Z $ and 
		$
		P_{\pm n} := \frac{\gamma}{2} \frac{M_n}{n} \pm M_n^{-1}$.
		These functions are the
   linear superposition  of  plane waves, traveling 
either to the right or to the left, 
see \cref{fig:StandingW} right. 
The condition $ \bar n_1 < \ldots < \bar n_\nu $ avoids 
the interference of plane waves with the same wave vector and amplitude, 
traveling in opposite directions, as 
$  \cos ( \bar n_a x + \Omega ( \bar n_a)  t) +
  \cos ( \bar n_a x - \Omega (- \bar n_a ) t) 	$, 
avoiding the formation of standing waves which are even in $ x $  at any  time $ t $,
as happens for $ \gamma =0 $.

\begin{theorem} \label{thm:main0}  {\bf (Quasi-periodic Stokes  waves \cite{BFM2})}
	Consider finitely many ``tangential sites'' 
	$  \bar n_1, \ldots, \bar n_\nu \subset \N $ satisfying 
	$ 1 \leq \bar n_1 < \ldots < \bar n_\nu $ 
	and signs $ \{ \sigma_1 , \ldots , \sigma_\nu \} $. Fix a subset 
	$ [\gamma_1, \gamma_2 ] \subset \R $.  Then 
	there exist $ \bar s >  0 $,  
	$ \varepsilon_0 \in (0,1) $ such that,  
	for any $ |\xi |   \leq \varepsilon_0^2  $, 
	$  \xi := (\xi_{ \sigma_a {\bar n}_a} )_{a = 1, \ldots, \nu} \in \R_+^\nu $, the following hold:
	\begin{enumerate}
	\item 
	there exists 
	a Borel set  $ {\mathcal G}_\xi \subset [\gamma_1, \gamma_2] 	$ 
	with density one at 
	$ \xi = 0 $, i.e. 
	$ \lim_{\xi \to 0} | {\mathcal G}_\xi |  = {\gamma}_2- {\gamma}_1 $;
	\item 
	for any  $ \gamma \in {\mathcal G}_\xi $,  the 
	gravity water waves equations \cref{eq:etapsi} 
	(with zero surface tension $ \kappa = 0 $) 
	have a  	quasi-periodic traveling wave solution of the form 
	\begin{equation}
		\label{QP:soluz}
			\begin{pmatrix}
				\eta ( t ,x) \\ \psi( t ,x)
			\end{pmatrix}  = \begin{pmatrix}
				\breve \eta_L ( \vec \omega t  - \vec \jmath x) 
				 \\ \breve \psi_L ( \vec \omega t  - \vec \jmath x)
			\end{pmatrix} + \breve r( \vec \omega t - \vec \jmath  x)   
	\end{equation}
	where $ (\breve \eta_L , \breve \psi_L) ( \vartheta ) $  are  defined in  \cref{sel.sol},  
	the wave vector is 
	 $ \vec \jmath  := (\sigma_a \bar n_a )_{a=1, \ldots,\nu} $, 
	\begin{enumerate}
	\item  the remainder 
	  $\breve r \in H^{\bar s} ( \T^\nu , \R^2)$ 	satisfies  
	$  \lim_{\xi \to 0} \frac{\| \breve r \|_{H^{\bar s}}}{\sqrt{|\xi|}} = 0 $; 
	\item 
	the frequency vector 
	$ \vec  \omega  := ( \omega_a)_{a=1, \ldots, \nu} 
	$ depends on $ \xi $, 	is Diophantine and satisfies
		$ \lim_{\xi \to 0}{\vec \omega} = 
	\vec \Omega  $. 
	\end{enumerate}
	\end{enumerate}
	In addition the quasi-periodic solutions  \cref{QP:soluz} are  linearly stable.  
\end{theorem}

An analogous KAM result with surface tension is proved in \cite{BFM1}. 
Let us make a comment on parameters. 
\\[1mm] 
{\sc Weak Birkhoff normal form.}
In \cref{thm:main0}, a small set of vorticities $ \gamma $ 
 is eliminated to impose Melnikov non-resonance conditions on the normal mode frequencies, cfr.
 \cref{dioph,1nd melnikov perd,2nd melnikov perd}. 
 Alternatively, one could eliminate  
 ``initial data'', 
 as inspired by the normal form approach for semilinear NLS in \cite{KP}.
For quasi-linear PDEs, this approach needs deep modifications, as 
  accomplished in  \cite{BBM16} for the KdV equation introducing 
  the concept of  {\it weak Birkhoff normal form}. This method was  successfully implemented 
  for deep water irrotational  waves  
  in \cite{FG}.

\subsection{Ideas of proofs.} \label{sec:1.3} The quasi-periodic solutions
of  \cref{thm:main00,thm:main0} are  constructed  by 
a Nash-Moser iterative scheme. The key analysis  concerns  
the linearized water waves equations \cref{eq:etapsi}  at 
an approximate quasi-periodic solution $ (\eta, \psi) ( \vec \omega t, x  )$. This is the 
time quasi-periodic linear  Hamiltonian system   
\be\label{LINOP}
\pa_t u + \begin{pmatrix} 
\partial_x \circ V + G(\eta) \circ B & - G(\eta) \\
(1 + B V_x) + B \circ G(\eta)  \circ B \  &  V \circ \partial_x - B \circ G(\eta) \end{pmatrix} 
u = 0
\ee
where $ (V,B) = (\Phi_x, \Phi_y)  $ is the velocity field 
evaluated at the free surface  $ (x,  \eta (\vec \omega t, x)) $, having 
small size $  \cO( \sqrt{|\xi |} ) $. 
Note that  the first order  transport operator 
$ \left(\begin{smallmatrix}
	\pa_x \circ V & 0 \\ 0 & V  \circ \pa_x
\end{smallmatrix}\right) $
in \cref{LINOP} 
is a singular perturbation of the unperturbed vector field 
$
\left( \begin{smallmatrix} 
0 & - G(0) \\
1  &  0 \end{smallmatrix}  \right) $, as we mentioned below \cref{diagNL0}. 

The key
result is  the {\it reducibility} of \cref{LINOP}   
   to a diagonal, constant coefficient linear system as \cref{eq:linKAM}.
This is achieved in two major conceptually different steps.  
For definiteness we consider   \cref{thm:main00}.

\begin{itemize}
\item[(a)]  {\bf (Reduction to constant coefficients up to smoothing remainders)}
As first task we conjugate \cref{LINOP}, via a  
quasi-periodically time-dependent change of variables acting in phase space  
$ v= \Psi(\vec \omega t ) u $ to a system of the form 
\be\label{op:redu1}
\partial_t v  \mathbb + \Big( \ii {\mathtt m}_{\frac12} |D|^{\frac12} \tanh^{\frac12}(\h |D|) + 
 \ii r (D) +    {\cal R}_0 (\vec \omega t) \Big) v = 0  
\ee
where $ {\mathtt m}_{\frac12} \sim 1 $ is a real constant  
independent of $ (t, x) $,  $ r (D) $ is real Fourier multiplier  operator 
independent of $ t$ of order $ -1/2 $,   
and the remainder $  {\cal R}_0  (\vec \omega t ) $ 
is a small quasi-periodic in time operator   
which is regularizing in space. 
This {\it pseudo-differential KAM normal form} is accomplished by an iterative process.
Let us outline only 
the 
straightening of the  highest order operator in \cref{LINOP},  namely
 the first order quasi-periodic transport 
Hamiltonian operator 
$
\pa_t + \left(\begin{smallmatrix}
	\pa_x \circ V & 0 \\ 0 & V  \circ \pa_x
\end{smallmatrix}\right) $. For standing waves, 
 $ V $ is odd in $ t$ and  in $ x $, and 
for $ V $  sufficiently small, 
this operator is ``rectified'' to  $ \pa_t  $
by conjugating with a family of circle diffeomorphism   
$ x \mapsto x + \beta (\vec \omega t, x )$ quasi-periodic in time. 
This is a small divisor problem, 
equivalent to proving that {\it all}  the solutions of the 
scalar characteristic ODE 
$ \dot x = V(\vec \om t, x )$ 
are quasi-periodic in time with frequency $ \vec \omega $, thereby 
forming a foliation of the torus $ \T^\nu $. 
In the traveling case, the operator  $
\pa_t + \left(\begin{smallmatrix}
	\pa_x \circ V & 0 \\ 0 & V  \circ \pa_x
\end{smallmatrix}\right)$ is rectified into $ \pa_t  + \mathtt m_1 \pa_x $ 
where $ \mathtt m_1 $ is a small constant, as the average of  $ V $. 
More generally we introduced  new quantitative {\it Egorov} theorems
for the reduction of pseudo-differential operators to constant  
Fourier multipliers, modulo smoothing remainders. 
These results show that a pseudo-differential operator, when conjugated under 
the flow of a first order transport-like  linear PDE, 
remains  pseudo-differential with a symbol 
having a decreasing-order expansion, up to a smoothing operator
that satisfies   precise  ``modulo-tame'' estimates. 
The most refined version is currently proved  in \cite{BFPT}. 
This effective use of 
Fourier integral operators (realized as flows of a transport PDE) is a 
crucial novel idea  for the development of 
KAM theory for quasilinear PDEs.
\item[(b)]  {\bf (KAM diagonalization in size)}
After \cref{op:redu1} is achieved, 
we implement  a KAM iterative scheme to reduce the perturbation
$ {\cal R}_0 ( \vec \omega t )$ 
 in size, 
yielding the linear autonomous diagonal system  \cref{eq:linKAM}. 
This is possible thanks to the Melnikov non-resonance conditions 
\cref{2nd melnikov perd} and because
${\cal R}_0 (\vec \omega t )  $ is sufficiently 
regularizing in space to compensate for the induced  loss of $ \cO(\perd)  $-derivatives.  
\end{itemize}

All the preceding changes of variables are bounded and satisfy ``tame'' estimates in Sobolev spaces. This allows to transfer tame estimates for the inverse 
of the diagonal autonomous linear system \cref{eq:linKAM}, into tame estimates 
for the inverse of the operator \cref{LINOP}  in the original physical coordinates,
sufficient for the convergence of a Nash-Moser iteration.

\smallskip 
These ideas and 
techniques provide a general set of tools for KAM theory for 
$ 1d$ quasi-linear  PDEs.   
For instance, 
they have been applied to quasi-linear perturbations of the Degasperis-Procesi  \cite{FGP1} 
and Klein-Gordon equations \cite{BFPT}; 
for 
the construction of large quasi-periodic multi-solitons for quasi-linear perturbations of 
KdV   \cite{BKM}; and to find 
quasi-periodic solutions of forced 3D Euler equations \cite{BM3d}.
For 2D Euler these techniques have been applied
to construct   
quasi-periodic vortex patches  
  near uniformly rotating Kirchhoff elliptical vortices   in \cite{BerHMasmo};
  vortex patches that desingularize point vortices  in \cite{HHM,HHR};  
and  stationary space-quasi-periodic solutions 
close to a Couette flow in \cite{FMM}. 
We quote \cite{HHM21,IGSP} 
 for KAM results for  $\alpha$-SQG equations.

\subsection{Open problems.} 
We list a series of open important  problems regarding KAM for water waves. 

\begin{enumerate}
\item[$P_{1,1}$)] {\it 3D fluids.} 
Traveling bi-periodic Stokes waves  were constructed 
in the gravity-capillary case in \cite{CN}  and in the pure gravity case
in \cite{IP-Mem-2009,IP2}. 
Quasi-periodic pure gravity traveling waves have been very recently proved 
in \cite{FMT}.   
The gravity-capillary case is still not understood. 
\item[$P_{1,2}$)] {\it Space quasi-periodic waves.} 
Are there time-periodic or quasi-periodic solutions to the water wave equations 
that are quasi-periodic in space? (not just periodic).

\item[$P_{1,3}$)] {\it Time Almost periodic solutions.} 
The water waves solutions in \cref{thm:main00,thm:main0} 
are time quasi-periodic. A very natural question regards 
the existence of more general almost-periodic solutions.  
\end{enumerate}

Actually results concerning  these problems would be very 
interesting also for any simpler PDE model.

\section{Almost global space periodic water waves.}\label{sec:BNF}

In the previous section we discussed the existence of many global in time 
space periodic solutions.   What can we say about the  time evolution of the 
other initial data? 
The  local   well posedness theory 
for free boundary Euler equations has been developed
in the pioneering works 
 \cite{Wu0,Wu1,Lannes, ABZ1,Lind1,Com,SZ1,IFT},  
proving the existence, for sufficiently regular initial data, of 
 classical  solutions on a small time interval.
 When specialized to initial data of size $ \e $  
 they  imply  a time of existence  larger than $ c \e^{-1} $. 
The natural question that arises in this context is the following:
\begin{itemize}
\item {\sc Question:}
Do space periodic water waves solutions exist for larger time intervals? 
\end{itemize}
 
The most general 
 almost global existence
 result 
  known so far is the following. Let $ \N_0 := \N \cup \{0\}$. 

\begin{theorem}\label{teo1}
{\bf (Almost global  in time 
gravity-capillary 
  water waves with constant vorticity \cite{BMM1})}
For any value of the gravity $ g > 0 $, depth $ \tth \in (0,+\infty] $ and vorticity $ \gamma \in \R $, there is a zero measure set $ {\mathcal K} \subset  (0,+\infty)$ such that, for any 
 surface tension coefficient $ {\kappa} \in (0,+\infty) \setminus {\mathcal{K}}  $, for any
$N$ in $\N_0 $, there is $s_0>0$ and, for any $s\geq s_0$, there are
$\e_0>0, c>0, C>0$ such that, for any $ 0 < \e < \e_0 $, any initial datum     
$$
(\eta_0,\psi_0) \in H_0^{s+\frac{1}{4}}(\mathbb{T},\R)\times {\dot H}^{s-\frac{1}{4}}(\mathbb{T},\R) \quad \mbox{ with} \quad
\| \eta_0 \|_{H^{s+\frac14}_0}+  \|\psi_0\|_{{\dot H}^{s-\frac14}}<\e \, , 
$$
system \cref{eq:etapsi} has a unique classical solution $(\eta,\psi)$ in 
\be\label{timeexi}
C^0\Bigl([-T_\e,T_\e],H_0^{s+\frac14}(\mathbb{T},\R)\times {\dot H}^{s-\frac14}(\mathbb{T},\R)\Bigr) 
\quad \mbox{with} \quad T_\e \geq c\e^{-N-1} \, , 
\ee
satisfying the initial condition $\eta |_{t=0} = \eta_0, \psi|_{t=0} = \psi_0$. Moreover
\be\label{boundsolN}
\sup_{t\in [-T_{\e},T_{\e}]}\big( 
 \|\eta\|_{H_0^{s+\frac14}}+\|\psi\|_{\dot{H}^{s-\frac14}}\big)\leq C\e\,.
\ee
\end{theorem}

Let us make some comments on the result. 
\\[1mm]
1. {\sc Standing waves.}
 This 
  result 
 is proved in Berti-Delort \cite{BD} for irrotational flows 
and  initial data $(\eta_0, \psi_0) $   
even in $ x $, so that  the 
solution in \cref{timeexi} remains even in $ x $ for all times (standing wave). 
In this subspace a reversible normal form 
is enough: reversibility ensures 
that the coefficient $ a $ of  resonant  monomial vector field  
$ a |u_{m_1}|^2 \ldots | u_{m_k} |^2 u_j $  is purely imaginary, making  
the actions $|u_j|^2 $  prime integrals of the normal form. 
In contrast,  
for proving the full  \cref{teo1},   it  is 
 necessary to 
 preserve  the Hamiltonian structure, 
 up to the degree  of homogeneity $ N $. 
\\[1mm]
2. {\sc Energy estimates.}
The  life span estimate $ T_\e \geq c \e^{-N-1} $ 
  and the bound \cref{boundsolN}
 for the solutions  follow by  
an energy estimate  for
$ \| ( \eta, \psi ) \|_{X^s} := \| \eta \|_{H_0^{s+\frac14}} +
\|\psi\|_{\dot{H}^{s-\frac14}} $  of the form 
\be\label{energyesint}
{\| (\eta,\psi)(t) \|}_{X^{s}}^{2} 
\lesssim_{s,N} {\|(\eta,\psi)(0)\|}^{2}_{X^{s}} +
 \int_{0}^{t} \| (\eta,\psi)(\tau)\|^{N+3}_{X^{s}} \, \di \tau 
\, . 
\ee
The fact that the  right hand side in \cref{energyesint} contains 
the same norm $ \| \ \|_{X^s} $ of the left hand side 
is  not trivial because  the equations \cref{eq:etapsi} are quasi-linear. 
The presence of the exponent $ N $   is  not trivial  because the  nonlinearity in 
\cref{eq:etapsi} vanishes only  quadratically for $ ( \eta, \psi ) = (0,0) $.
It is the outcome of the paradifferential Hamiltonian  BNF. 
\\[1mm]
3.  {\sc Non-resonance conditions.}
 The restriction on the  surface tension  required in \cref{teo1} 
arises  to ensure the absence of $N$-wave resonant interactions  among the 
linear frequencies $ \Omega (j) $  in \cref{omegonejin} 
\be\label{Nwavei}
\big| \Omega (j_1) + \ldots +\Omega (j_p)-\Omega (j_{p+1})-
\ldots-  \Omega (j_N) \big| \gtrsim 
 \max(|j_1|, \ldots, |j_N|)^{-\tau} \, , 
\ee
for any $ 1 \leq p \leq N $,  and 
wave vectors  $ j_1 , \ldots 
j_N  $ satisfying 
$ \{ |j_1|, \ldots, | j_p|\}\not=\{ |j_{p+1}|, \ldots, |j_N|\} $. 
Note that 
 if $ N $ is even and $ p = N / 2 $  then the left hand side in \cref{Nwavei} may 
 vanish for 
$ \{ |j_1|, \ldots, | j_{\frac{N}{2}}|\} =\{ |j_{\frac{N}{2}+1}|, \ldots, |j_N|\} $
(this is actually the case if $ j \mapsto \Omega (j) $ is even, i.e. for irrotational waves). 
Because of these exact resonances some 
``resonant'' vector field can not be eliminated in the normal form.
The small divisors  conditions \cref{Nwavei} 
 are very weak, as the  lower bound  
 depends on the highest wave vector,  
causing a loss of  $\tau $-derivatives in  the normal form 
 process.
This is compensated 
by a  paradifferential reformulation of the  equations and a 
regularization procedure that transforms them 
 into constant coefficient 
symbols, up to smoothing remainders, as we  comment  below.
\\[1mm]
\noindent 
4. {\sc Long time existence of space periodic 
water waves.} We now describe other 
long time existence  results  
for space periodic water waves proved in the literature.  
\\[1mm]
($i$)  $ T_\e \geq c \e^{-2} $. 
The works  \cite{Wu0,IP,AlDe1} 
for pure gravity waves, and  \cite{ifrTat,IP3} for 
  $\kappa>0$, $g=0$ and $ \tth = + \infty $,  prove that  
small data
of size $ \e $ 
(periodic or on the line)
give rise to  irrotational water waves 
defined on a time interval   
at least $c \e^{-2}$. 
We quote  
\cite{IFT} 
for  $\kappa = 0$, $g>0$, $ \tth = + \infty $
and constant
vorticity,  \cite{HuIT} for irrotational fluids, and 
\cite{HarIT}  in finite depth. 
All the  previous results hold in absence of three wave interactions.  
Exploiting the Hamiltonian nature of the  equations,  \cite{BFF}  
proved, for any value of gravity, surface tension and depth,  
that  $ T_\varepsilon \gtrsim  \e^{-2}$. 
Interestingly in  these cases non trivial, finitely many, 
three wave interactions 
$ \Omega (j_1)  \pm \Omega (j_2)\pm
\Omega (j_3) = 0  $, 
$ j_1 \pm j_2 \pm j_3  = 0 $, 
may occur,  that in fluid mechanics give rise to the well known Wilton ripples instabilities. 
\\[1mm]($ii$)   $T_\e \geq c \e^{-3}$.
The pure 
gravity water waves equations in deep water have been proved 
to be integrable  up to 
quartic vector fields in  \cite{BFP}, 
 implying in particular the well-posedness of the solutions 
for times $ \gtrsim  \e^{-3} $.  The  approach in  \cite{BFP}   
 to  rigorously   justify the formal computations in  \cite{ZakD} 
is  based on a uniqueness argument of the normal form: 
 the paradifferential BNF and the formal one 
  have to coincide up to cubic degree.
 The result has been recently extended in \cite{Wu3} through a novel configuration space approach, applicable to a broader class of initial data, and in 
\cite{DIP}  for waves with large period or  initial data in $ H^s (\R) $, obtaining, in this case, 
a time of existence $ \gtrsim \e^{-6} $.

\subsection{Ideas of proofs.} \label{sec:ideasBNF}
The  {\it loss of derivatives}  caused, along the normal form procedure,  
by the quasi-linearity of the equations and the small divisors induced by 
the $ N$-wave nearly resonances \cref{Nwavei},  
is addressed by the idea of  {\it paradifferential Birkhoff normal form}  introduced in \cite{BD}, 
that we now describe. The first step is to ``paralinearize'' 
the Hamiltonian PDE \cref{eq:etapsi} that,  in the complex 
variable  \cref{defMM-1}, may be written  as 
 \be\label{JcH.intro}
\pa_t U = {\cal J}  \nabla H (U)  \qquad \text{where} \qquad  {\cal J} := \begin{pmatrix}0 & - \im \\ \im  & 0  \end{pmatrix}  \, , 
\quad U = \begin{pmatrix}u \\ \bar u \end{pmatrix} \, ,  
\ee
and  $H $ is the Hamiltonian \cref{H.gamma}. 
This means to write \cref{JcH.intro} in the form  
\be\label{equaz1}
\pa_t U =   {\cal J} \Opbw{A(U; x,\xi)}[U] + R(U)[U] 
\ee
where $ A(U;x,\xi)$ is a matrix of symbols and $ R(U) $ are $ \varrho$-smoothing operators.
The symbols $ A (U; x, \xi )$ depend  on the variable $ U $ which has
 only finite space regularity. 
All the symbols and smoothing remainders admit a Taylor expansion in homogenous components in 
$ U $   up to degree $ N $. Their precise definitions are rather technical, but we say the following.  A homogenous symbol 
$ a_q(U;x,\xi)$ of degree $ q \in \N $ and  order $ m\in \mathbb{R} $, has the form 
\be\label{sviFou}
a_q(U;x,\xi)=  \!\!\!\!
\sum_{j_1, \ldots, j_q \in {\Z \setminus
 \{ 0 \}},\sigma_1, \ldots, \sigma_q\in \{\pm 1\}}  \! \! \! \! \!\!\!\!\!  
 \left( a_q\right)_{j_1,\ldots,j_q}^{\sigma_1,\ldots, \sigma_q}(\xi) 
u_{j_1}^{\sigma_1}\dots u_{j_q}^{\sigma_q}
 e^{\im (\sigma_1 j_1 + \ldots + \sigma_q j_q) x}  
\ee
where 
$ \left( a_q\right)_{j_1,\ldots,j_q}^{\sigma_1,\ldots, \sigma_q}(\xi)  
$ are   scalar valued  Fourier multiplier 
of order $m$.
The Bony-Weyl {\it paradifferential} operator is defined as the Weyl quantization of a regularized 
$ C^\infty $ symbol
$a_\chi (x, \xi)  $,    acting on a $ 2 \pi $-periodic function $u(x)$ as 
$$
\left({\rm Op}^{BW} {(a)}[u] \right)(x) 
:= {\rm Op}^W {(a_\chi )}[u]  := {\mathop \sum}_{j \in \Z} \Big( {\mathop \sum}_{k \in \Z}
\widehat a_\chi \big(j-k, \tfrac{k + j}{2}\big)  \, u_k \Big) e^{\im j x } \, ,
$$
where $ \widehat a_\chi (n, \xi) $ is the $n$-th Fourier coefficient of the $ 2 \pi $-periodic function 
$ x \mapsto a_\chi (x, \xi)  $, and  $ a_\chi (x, \xi) $
is achieved restricting the symbol in \cref{sviFou}
to  Fourier harmonics with {\it low} frequencies 
 $ | j_1 |, \ldots, | j_q|  \leq \delta \langle \xi \rangle $ where  $ \langle \xi \rangle  := \max (1, |\xi| ) $. 
If  $ \delta > 0 $ is sufficiently small, there is 
a {\it symbolic calculus} for paradifferential operators  analogous to that of pseudo-differential operators. 
A $ \varrho $-smoothing operator of homogeneity $ q $ satisfies 
estimates as 
$$
 \| R(U)V \|_{\dot{H}^{s + \varrho}} 
   \lesssim_{s,q}  \|{V}\|_{\dot{H}^{s}}\|{U}\|_{\dot{H}^{s_0}}^{q+1} 
 +\|{V}\|_{\dot{H}^{s_0}}\|U\|_{\dot{H}^{s_0}}^{q}\|{U}\|_{\dot{H}^{s}} \, , \quad \forall s \geq s_0 := s_0 (\varrho) \, . 
$$
The paradifferential Hamiltonian BNF 
is achieved in three major conceptually different steps.  
\begin{itemize}
\item[(a)] {\bf (Para-differential reduction to constant coefficients symbols up to smoothing remainders)}
We perform several transformations of the form 
\be\label{WphiU}
W = \Phi (U) U  \, , \quad  \text{where} \quad \Phi (U) : \dot H^s (\T, \C^2) \to \dot H^s (\T, \C^2) \, , 
\ \forall s \in \R  \, , 
\ee
(including the Alinhac good unknown \cite{Lannes,ABZ1})
to convert \cref{equaz1} into an equation  
\be\label{st1nfr}
\pa_t w = \Opbw{
 \ii \tm_{\frac32} (U;\xi)}w + R(U)U \, , \quad W = 
 \big( \begin{smallmatrix}w \\ \bar w \end{smallmatrix}\big)  \, , 
\ee 
where  
 $ \tm_{\frac32}(U,\xi) $ is  a real symbol of order $ 3/ 2 $, independent of $ x $, 
of  the form 
$$
 \tm_{\frac32}(U;\xi):= -  (1+\zeta(U))
 \sqrt{ \xi \tanh(\tth \xi) \big(g+ \kappa \xi^2 + \tfrac{\gamma^2}{4} \tfrac{ \tanh (\tth \xi)}{\xi} \big)} - \tfrac{\gamma}{2} \tanh (\tth \xi) -  \mathtt V(U) \xi + \text{lower order symbols} \, ,
$$
with real functions $\zeta(U)$, $\mathtt V(U)$   independent of $ x$, 
vanishing at $ U = 0 $. Thus, for $ U = 0 $, 
the symbol $   \tm_{\frac32}(U;\xi) $  reduces   
to the unperturbed dispersion relation $ - \Omega (\xi)  $ in 
\cref{omegonejin}. 
The multiplicative 
constant $ \zeta (U)   $
takes into account the quasi-linear  effects of the nonlinearity. 
The reduction of the highest-order term is 
obtained defining a ``para-composition'' operator
as  the  flow of a first-order para-transport PDE. 
In  \cite{BD}  we prove
that a para-differential operator, when conjugated under such a flow, 
remains  para-differential,  with a symbol 
having an expansion in  decreasing orders.   
This is the paradifferential analogue of the Egorov theorem.

By ignoring the remainder term 
$R(U)U $  in  \cref{st1nfr}, the $ L^2 (\T) $ and the $ H^s (\T) $  norms 
of the solutions  are preserved indefinitely, as the symbol $ \tm_{\frac32}  (U; \xi) $ is 
real, so $ \Opbw{
 \ii \tm_{\frac32} (U;\xi)}$ is self-adjoint,  
 and independent of $ x $. However, the quadratic term 
$ R(U) U $  still  limits the stability  time to  $ \e^{-1}$, so we must  eliminate it. 
\item[(b)] {\bf (Symplectic correctors up to homogeneity $ N $)} 
Substituting in  \cref{st1nfr} the relation  $ U = \Phi^{-1} (W) W + \ldots $ 
obtained  inverting \cref{WphiU} 
 approximately up to homogeneity $ N $, we obtain a similar  equation 
\be\label{st1nfrnew}
\pa_t w = \Opbw{\ii \tilde \tm_{\frac32}(W,\xi)} w 
 + \tilde R(W) W  \, ,
\ee
which 
 is not Hamiltonian because \cref{WphiU} is not symplectic. A major contribution of \cite{BMM1} 
is to restore the Hamiltonian nature of the equation \cref{st1nfrnew} up to homogeneity $ N $, 
by composing \cref{WphiU}  with symplectic correctors, 
which are 
smoothing perturbations of the identity. The symplectic correctors
are constructed by a novel perturbative procedure inspired by the 
celebrated Moser proof of the Darboux theorem. In this way we obtain 
an equation as \cref{st1nfrnew} which is Hamiltonian  up to degree of homogeneity $ N $. 
\item[(c)]  {\bf (Paradifferential Hamiltonian BNF)} 
Thanks to the non-resonance conditions \cref{Nwavei}  we eliminate 
the Hamiltonian monomials up to homogeneity $ N $ --both in the symbols and in the  smoothing remainder in \cref{st1nfrnew}--  
that do not Poisson commute with the ``super-actions'' 
$   |z_{n}|^2 + |z_{-n}|^2  $ for any $ n \in \N $. This is valid for positive surface tension.
These changes of variables  are bounded in a sufficiently smooth Sobolev space.  
The remaining monomials  in the resonant normal form allow
an exchange of energy between the Fourier modes $ \{ z_{n}, $ $ z_{-n} \} $.
However, thanks to  its Hamiltonian structure, the super-actions 
 are prime integrals of the normal form, and  the Sobolev norms of the solutions are preserved.
 This BNF procedure 
 ultimately implies an energy estimate as \cref{energyesint}. 
\end{itemize}

This set of  ideas   and  techniques provide a general framework  for normal form theory for 
$ 1d$ quasi-linear Hamiltonian PDEs.  
For instance they have been 
 applied to prove a  quadratic lifespan result
 \cite{BCGS} for vortex patches of $ \alpha $-SQG near  circles
 and in \cite{MRS}  almost global existence of 
 Kelvin-Helmholtz vortex sheets near a circle.

\subsection{Open problems.} 
We list a series of open problems that we think are very interesting. 

\begin{enumerate}
\item[$P_{2,1}$)] {\it Global existence.}
Are the 
solutions proved in  \cref{teo1} 
global in time or not? 
Global well-posedness for  \cref{eq:etapsi} 
is an open problem, even for smooth and  localized irrotational initial data in $ \R^2 $ 
 (\cite{DIPP} applies in  $\R^3 $).
\item[$P_{2,2}$)] {\it 3D fluids.} 
In higher space dimension 
 the paradifferential reduction up to smoothing remainders is no longer available. 
We quote the $ \e^{- \frac53+} $ long time existence result in 
 \cite{IoP} for pure gravity periodic waves with $ x \in \T^2 $.  
\item[$P_{2,3}$)] {\it Pure gravity water waves.} In \cref{teo1} we exploit 
capillarity in order to prove the reduction to constant coefficient symbols.  
For pure gravity waves we do not have an almost global existence result. 
\item[$P_{2,4}$)] {\it Rational normal form.} In the deep water pure gravity case \cite{BFP} the quartic
Hamiltonian is integrable and it is tempting  to prove the existence of solutions
for longer times $ \gtrsim \e^{-4}, \e^{-5}, \ldots $ erasing initial conditions, in the spirit of the 
rational normal form result \cite{BFG} for semilinear NLS.   
\end{enumerate}

\section{Modulational instability.}\label{sec:BF}
The Stokes waves \cref{eq:Stokes1}-\cref{eq:Stokes2} 
 are often used in oceanography to model 
 complex nonlinear   wave interactions,  and 
a question of fundamental physical importance concerns their stability or instability,  
when subjected to wave disturbances of different spatial periods. 

To mathematically describe this phenomenon, we   linearize the  water waves equations at the Stokes waves in the reference frame moving with the speed $ c_\e $ of the wave.   
After some changes of variables (including a linear Alinhac 
good unknown and a conformal mapping),
this is  the  linear autonomous real Hamiltonian system  
\begin{equation}
\label{cLepsilon}
\pa_t h (t,x) = \mathcal{L}_{\e}
 \, [h(t,x)]    \qquad 
\text{where}  \qquad 
\cL_\e := 
\begin{pmatrix} \pa_x \circ (\tc_\tth+p_\e(x)) &  |D|\tanh((\tth+\mathtt{f}_\e) |D|)  
 \\ - (1+a_\e(x)) &   (\tc_\tth+p_\e(x))\pa_x \end{pmatrix} \, , 
\end{equation}
$ \ch $ is the linear  speed   in \cref{eq:Stokes2}, 
$\mathtt{f}_\e = \cO(\e^2 ) $ is a real number 
analytic in $ \e $,  whereas 
$p_\e (x), a_\e (x) $ are even real analytic $2\pi$-periodic  functions of size $ \e $
analytic in $\e$ as well (they depend on the Stokes wave). In  deep water $ \tth = + \infty $ just replace 
$ \tanh((\tth+\mathtt{f}_\e) |D|)  $ with $ 1 $. 

The dynamics of \cref{cLepsilon} for    
perturbations 
$ h (t,\cdot) $  in  $ L^2 (\R) $, resp. in 
 $  L^2 (\T_N)$ 
where $ \T_N := \R / (2 \pi N \Z ) $ is 
a large torus, is 
fully determined by answering the following
\begin{itemize}\item
{\sc Question:} {\it What is the $ L^2 (\R)$, resp. $ L^2 (\T_N) $, spectrum of
$ \mathcal{L}_{\e}$?} 
\end{itemize}

\begin{figure}  \label{fig:spectrum}
\boxed{
\begin{minipage}{0.5\textwidth}
\vspace{0.47cm}
\centering
\includegraphics[width=6.6cm]{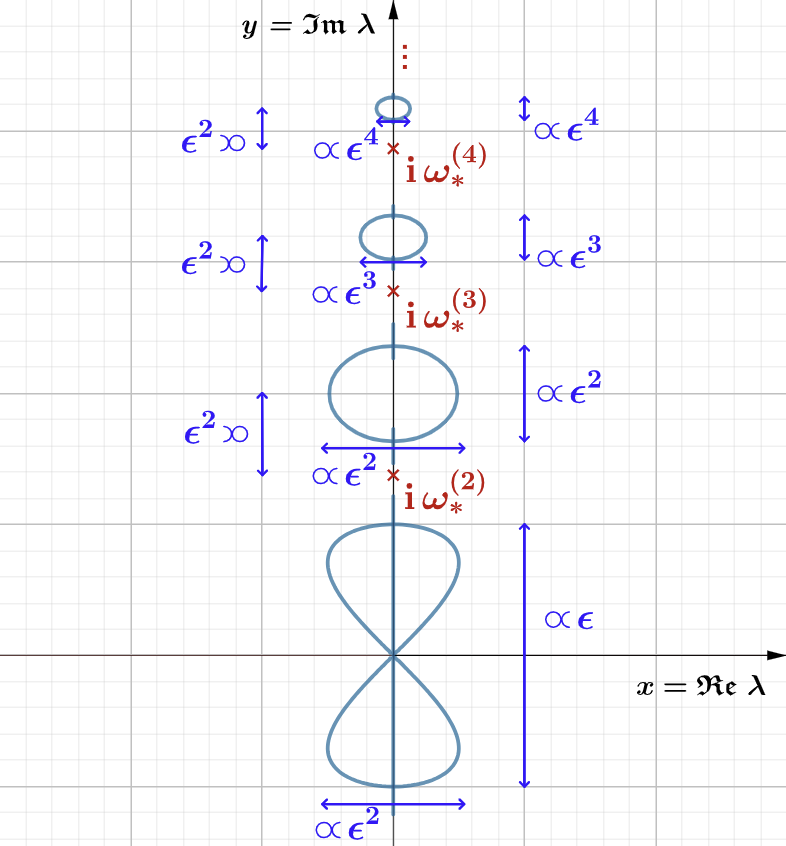}
\vspace{0.47cm}
\end{minipage}}
\boxed{\begin{minipage}{0.47\textwidth}
\centering
\includegraphics[width=4.35cm]{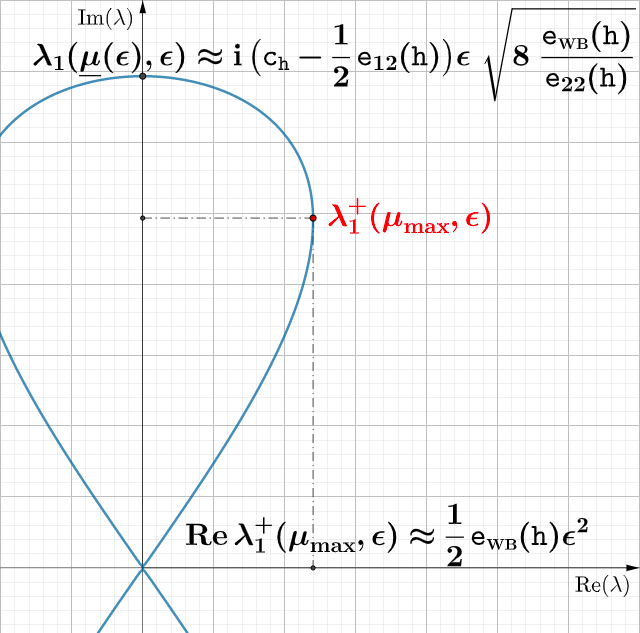}
\rule{\textwidth}{0.4pt}
\includegraphics[width=4.35cm]{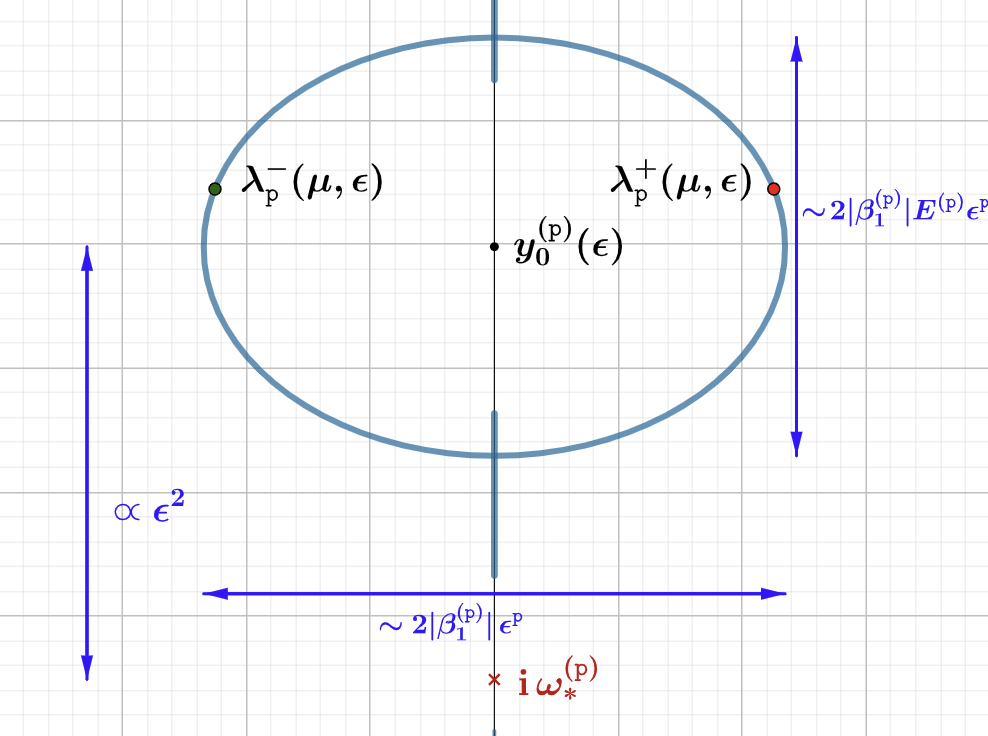}
\end{minipage}}
    \caption{
    spectral bands with non zero real part of the 
$ L^2 (\R) $-spectrum of
$ \cL_{\e} $. The ``figure 8'' is present for $ \tth > 1.363... $.
For any 
$ {\mathtt p} \geq 2 $,
for any depth except a closed set, there exists 
an isola of instability with size 
$ \propto  | \beta_1^{(\tp)}(\tth)| \e^{{\mathtt p} } $. Its
center $ y_0^{({\mathtt p} )} (\e) $ is $ \cO(\e^2) $ distant from $ \im \omega_*^{({\mathtt p} )}  $.
The notation $ a \propto b $ means that $ a b^{-1} \to c  $ as $ \e \to 0 $. 
}
\end{figure}

After preliminary studies 
conducted by  Philipps and McLean  in the '70-'80, 
the paper \cite{DO} 
numerically neatly revealed 
 two distinct types of unstable structures in the spectrum 
as shown in \cref{fig:spectrum}: 
a ``figure 8'' shaped band of unstable spectrum crossing the origin --corresponding to the
Benjamin-Feir instabilities-- and  
elliptical ``isolas'' of smaller 
size centered on the imaginary axis, giving rise to other high frequency  
instabilities.  
The first  rigorous  mathematical 
proofs of the Benjamin-Feir  
instability  
of the Stokes waves under  longitudinal perturbations 
 were obtained   in \cite{BrM}  in finite depth 
and, recently, in deep water,   in \cite{NS}. 
These results describe the local cross of the figure 8 near zero. 
The first isola closest to zero was described analytically in finite depth in  \cite{HY}.

In the last few years,  in a series of papers starting with \cite{BMV1}, we obtained: 
\begin{enumerate}
\item {\bf (Unstable spectrum near zero)}
the complete description of the $ L^2 (\R) $-spectrum of the linearized operator $ \cL_\e $ near the origin. This was done in deep water \cite{BMV1}, finite depth \cite{BMV3}, and, in  \cite{BMV_ed}, at the
famous Whitam-Benjamin critical depth 
$\tth_{{\mathtt W}{\mathtt B}} = 1.363... $,  below which the Stokes waves become linearly stable.
\item {\bf (Unstable spectrum away from zero)} The existence of  arbitrarily  
many high-frequency unstable spectral isolas away from the origin, for any depths
$ \tth > 0 $  except a closed set, see  \cite{BCMV}. These instabilities  become increasingly
stronger in the shallow water 
limit $  \tth \to 0^+ $. 
\end{enumerate}

These results 
provide a rigorous justification for the picture of  \cref{fig:spectrum} 
for finite depths. 

Modulational instability  under transversal wave perturbations is  proved in
\cite{CNS,JRSY,HTW} and references therein. 

\smallskip 

We first provide a precise account of the results  \cite{BMV1,BMV3,BMV_ed}.
Since the operator $\cL_\e $ in 
\cref{cLepsilon} has $2\pi$-periodic coefficients, its $ L^2(\R) $ spectrum is 
effectively described using  the Bloch-Floquet decomposition 
\be\label{eq:BF}
     \sigma_{L^2(\R)}\big( 
     \cL_\e  \big)
 = \bigcup_{\mu \in \big[-\tfrac12, \tfrac12 \big)}
 \sigma_{L^2(\T)}\big( \cL_{\mu,\e} 
 \big)\, \qquad \textup{where}\qquad 
 \cL_{\mu,\e} := 
e^{-\im \mu x}\, \cL_\e\,  e^{\im \mu x}  
\ee
is 
the complex  Hamiltonian and reversible pseudo-differential  operator
\begin{align}\label{WW}
 \cL_{\mu,\e}    
 = 
 \underbrace{\begin{pmatrix} 0 & \uno\\ -\uno & 0 \end{pmatrix}}_{=J} 
 \underbrace{\begin{pmatrix} 1+a_\e(x) & -(\ch+p_\e(x))(\pa_x+\im \mu) \\ (\pa_x+\im\mu)\circ (\ch+p_\e(x)) & |D+\mu| \tanh\big((\tth + \ttf_\e) |D+\mu| \big) 
 \end{pmatrix} }_{=: \cB_{\mu,\e}}  \, .   
\end{align} 
Here complex Hamiltonian means that 
$ \cB_{\mu,\e} = \cB_{\mu,\e}^*   $ 
where  $\cB_{\mu,\e}^*$ 
 is the adjoint with respect to the complex scalar product of $L^2(\T, \C^2)$ and 
 reversible means   that
$$
 \cL_{\mu,\e}\circ \bro =- \bro \circ \cL_{\mu,\e} \qquad \text{ where} \qquad 
   \bro \vet{\eta(x)}{\psi(x)} := \vet{\bar\eta(-x)}{-\bar\psi(-x)} \, . 
$$
The spectrum  of $ \cL_{\mu,\e} $  on $ L^2 (\T ) $ 
is discrete. As the {\it Floquet parameter}  $ \mu $ varies, 
the   eigenvalues  of $ \cL_{\mu,\e} $ trace out  the continuous {\it spectral  bands} of 
$ 	\sigma_{L^2(\R)} ( \cL_\e ) $. 
The  spectrum $ 	\sigma_{L^2(\T_N)} ( \cL_\e ) $ 
 is obtained in a similar manner, 
but the Floquet parameter $ \mu $ 
is restricted to a finite set of values separated by  gaps of size  approximately 
$ N^{-1} $. 
As the size of the  large box torus   $ N \to + \infty $ the eigenvalues
of $ \cL_{\mu,\e} $  become increasingly dense.

Any solution  of \cref{cLepsilon} can be decomposed into a linear superposition of Bloch-waves
\begin{equation}\label{hesplode}
h(t,x) = e^{ \lambda t} e^{\im \mu x} v(x)
\end{equation}
where  $\lambda $ is an eigenvalue of $\cL_{\mu,\e} $ 
with  corresponding   eigenvector  $ v(x) \in L^2(\T) $. 
 If the real part of $\lambda $ is  positive,  
 the solution \cref{hesplode}  grows exponentially fast in time,
 and   we say that $ \lambda $ is {\it unstable}.  
 
We also remark that the spectrum $ \sigma_{L^2(\T)}(\cL_{\mu,\e}) $
is a set which is $1 $-periodic in $ \mu $. This is why in \cref{eq:BF}
it is sufficient to consider $ \mu $ in  the ``first zone of Brillouin'' 
$ [-\tfrac12, \tfrac12)$, mod $ \Z $. 

\smallskip
Since $\cL_{\mu,\e}$ is Hamiltonian, eigenvalues with non zero real part may arise only from multiple
  eigenvalues of $\cL_{\mu,0}$, because   if $\lambda$ is an eigenvalue of $\cL_{\mu,\e}$ then  also $-\bar \lambda$ is.  
  \\[1mm]
  {\sc The spectrum of $\cL_{\mu,0}$.} 
  Let us determine the multiple eigenvalues of the Fourier multiplier matrix operator  
$$
 \cL_{\mu,0} := 
 \begin{pmatrix} 
 \ch (\pa_x+\im\mu)  & |D+\mu|\tanh(\tth |D+\mu|) \\ -1 & \ch(\pa_x+\im\mu) \end{pmatrix} \, .  
$$
The spectrum of  $\cL_{\mu,0}  $ on $L^2(\T,\C^2)$  consists of the   purely imaginary numbers
$$
\lambda_j^\sigma(\mu):=
\im \omega^\sigma(j+\mu,\tth) 
\quad \  \forall \sigma   = \pm \, , \,  j\in \Z \, , \quad
\text{where} \quad   
\omega^\sigma(\varphi,\tth) :=  \ch \varphi -\sigma \sqrt{\varphi\tanh(\tth \varphi)} \, , \  \varphi \in \R \, . 
$$
For $\mu=0$ the operator $\cL_{0,0}$ possesses the multiple 
zero  eigenvalue 
$ \lambda_0^+(0) = \lambda_0^-(0) = \lambda_1^+(0) = \lambda_{-1}^-(0)=0 $.

\begin{remark}
The unperturbed 
eigenvalue $ 0 $ remains an eigenvalue of $ \cL_{0,\e} $ with algebraic multiplicity $ 4 $, 
for any $ \e $, because Stokes waves are not isolated but 
rather appear on a $ 4 $-dimensional  manifold. For instance, 
in infinite depth,  they can be written as 
$ \big( \breve \eta_\varepsilon + \frac{P}{g}, \breve \psi_\varepsilon + \psi_0 \big) (x+ \theta) $, for arbitrary  real parameters 
$ \varepsilon, P  \, , \theta  $ and $ \psi_0 $. 
\end{remark}

For any  $ \mu \in \R$ and  $\tth > 0 $,  
the non zero eigenvalues $ \lambda_j^\sigma(\mu)  $ of  $\cL_{\mu,0} $ are 
either simple or  double. More precisely 
there exists a sequence of Floquet parameters 
\be\label{eq:umu}
 \umu=\uphi(\tp,\tth) > 0  \, , \mod \Z \, , \quad \tp = 2,3, \ldots\,  , 
 \ee
for which two eigenvalues  $  \lambda_j^\sigma(\mu)  $ of $\cL_{\umu,0} $ ``collide'' into 
a double eigenvalue  $ \im \omega_*^{(\tp)}(\tth ) $. These  double  eigenvalues  
form a diverging sequence  
\be\label{eq:omegastar}
\{\pm \im \omega_*^{(\tp)}(\tth ) \}_{\tp=2,3,\dots} \, ,
\quad 
0<\omega_*^{(2)}(\tth )  
< \dots < \omega_*^{(\tp)}(\tth ) < \dots \, , \qquad 
\lim_{\tp \to + \infty}
\omega_*^{(\tp)}(\tth) =  + \infty \, . 
\ee
\begin{itemize}
\item 
{\sc Question:} For small $ \e \neq 0 $, do unstable eigenvalues of $ \cL_{\mu,\e} $ emerge from 
 the multiple eigenvalue $ 0 $ of the  unperturbed operator $ \cL_{0,0} $, for  $ \mu $ near $0 $? And  
 from the   double  eigenvalues $ \{\pm \im \omega_*^{(\tp)}(\tth ) \}_{\tp=2,3,\ldots} $ 
 in \cref{eq:omegastar}
of $ \cL_{\umu,0}$, for  $ \mu $ near $ \umu  $ ?
 \end{itemize}

\subsection{Unstable spectral bands close to zero: the ``figure 8''.}
The next result describes the  unstable spectral bands in the spectrum 
of $ \cL_\e $ close to zero in  deep water 
$ \tth = + \infty $.  We denote by $r(\e^{m_1} \mu^{n_1}, \ldots, \e^{m_q} \mu^{n_q})$
a real analytic function fulfilling  for some $C >0$ and $\e, \mu$  small,  $| r(\e^{m_1} \mu^{n_1}, \ldots, \e^{m_q} \mu^{n_q}) | \leq 
C \sum_{j=1}^q
 |\e|^{m_j} |\mu|^{n_j}
$.
\begin{theorem}\label{TeoremoneFinale}
{\bf (Benjamin-Feir spectrum in deep water \cite{BMV1})}
There exist $ \e_0, \mu_0 > 0 $  and 
an analytic  function 
 $$ 
\underline \mu: [0,\e_0)\to [0,\mu_0) \, , \quad  \underline  \mu(\e)  =   2\sqrt{2} \e(1+r(\e)) \, , 
 $$ 
 such that,  for any  $  \e \in [0, \e_0)  $, the 
 operator  $\cL_{\mu,\e}$  
 has two eigenvalues
  $\lambda^\pm_1 (\mu,\e)$ of the form 
\begin{equation}\label{eigelemu}
 \begin{cases}
  \frac12\im\mu+\im r(\mu\e^2,\mu^2\e,\mu^3)\pm
  \frac{\mu}{8}\sqrt{8\e^2\big(1+r_0(\e,\mu)\big)-\mu^2\big(1+r_0'(\e,\mu)\big)} \, ,  & \forall   \mu \in [0, \underline \mu (\e)) \, , \\
 \frac12\im \underline \mu (\e)+\im r(\e^3) \, ,  & \mu= \underline \mu  (\e) \, , \\
 \frac12\im\mu+\im r(\mu\e^2,\mu^2\e,\mu^3)\pm\im\frac{\mu}{8}\sqrt{\mu^2\big(1+r_0'(\e,\mu)\big)-8\e^2\big(1+r_0(\e,\mu)\big)} \, ,  & \forall \mu \in ( \underline \mu  (\e), \mu_0) \, .
\end{cases}
\end{equation}
The function 
$$ 
8\e^2\big(1+r_0(\e,\mu)\big)-\mu^2\big(1+r_0'(\e,\mu) )  >0 \, , \quad 
respectively <0 \, , 
$$
 provided $0<\mu < \underline{\mu}(\e)$,  respectively $\mu > \underline{\mu}(\e)$.
\end{theorem}

For values of the Floquet parameter $ \mu \in (0, \underline \mu (\e)) $ the eigenvalues 
$\lambda^\pm_1 (\mu, \epsilon) $ in \cref{eigelemu} trace the upper part of a 
{\it figure 8}; 
  at $ \mu =  \underline  \mu (\e)$ they
 collide on the imaginary axis {\it far} from $ 0 $ and  they remain purely imaginary 
 for $ \mu > \underline  \mu (\e) $,
see \cref{fig:spectrum}. The maximum amplitude of the figure 8 
 is reached at a complex number with real part $\sim 
    \tfrac12  \e^2 $.  
The upper part of   the figure  8 is 
well-approximated by the curve 
$\mu \mapsto (\pm\frac{\mu}{8}\sqrt{8\e^2 -\mu^2}, \frac12 \mu)$. 
The higher order corrections of the eigenvalues $ \lambda_1^\pm (\mu,\e ) $ in \cref{eigelemu} 
are algorithmically  computable. For $ \mu < 0 $ the operator $ {\mathcal L}_{\mu,\e} $
possesses symmetric  eigenvalues 
in the lower half complex plane.

\smallskip

 The complete description of the Benjamin-Feir spectrum 
in finite depth   is given 
in  \cite{BMV3} and  \cite{BMV_ed}: 
\begin{itemize}
\item {\bf Sufficiently deep water \cite{BMV3}.} 
   If  $\tth >   \tth_{{\mathtt W}{\mathtt B}} = 1.363...$, 
  for any $ \e \in (0, \e_0)  $ small (depending on $ \tth $), the spectrum 
  $ \sigma_{L^2(\R)} (\cL_\e )$ contains
a ``figure 8''  region, 
 with maximum amplitude 
  $  \sim  \tfrac12  {\mathtt e}_{{\mathtt W}{\mathtt B}} (\tth) \e^2$
  (see \cref{fig:spectrum}, top-right), where the constant 
    $ {\mathtt e}_{{\mathtt W}{\mathtt B}} (\tth) \to  1^-  $ as $ \tth \to + \infty $
    and  $ {\mathtt e}_{{\mathtt W}{\mathtt B}} (\tth) \to  0^+ $ 
     as $ \tth \to \tth_{{\mathtt W}{\mathtt B}}^+ $.
\item {\bf Shallow water \cite{BMV3}.} 
If $ 0 < \tth <  \tth_{{\mathtt W}{\mathtt B}} $,  for any $ \e \in (0, \e_0) $ small 
(depending on $ \tth $) 
 the spectrum  $ \sigma_{L^2(\R)} (\cL_\e )$  
 near zero  is purely imaginary.
\item  
{\bf Critical depth  $\tth =  \tth_{{\mathtt W}{\mathtt B}} $ \cite{BMV_ed}.} 
  At the critical depth  $\tth =  \tth_{{\mathtt W}{\mathtt B}} = 1.363... $
 predicted by 
  Whitham and Benjamin (and for values of 
  $ (\tth,\e) $ satisfying
   $ \tth > \tth_{{\mathtt W}{\mathtt B}}  - c \e^2 + \cO(\e^3)$, $ c >0 $)  the Stokes waves become modulationally unstable. 
This  result solved a long-standing and controversial question 
within the fluid dynamics community.   
\end{itemize}

\subsection{Ideas of proofs.} \label{sec:ideasBF}
We 
 implement Kato's theory for separated eigenvalues for 
 the  complex Hamiltonian and reversible  
operators $\cL_{\mu,\e}$.
 The quadruple zero eigenvalue of $\cL_{0,0}$ is isolated, meaning there is  
 a  path  $\Gamma$ separating   $ 0 $ from the rest of the spectrum
 $ {\cal L}_{0,0} $. For  sufficiently small $ \e, \mu $,  the  operator
$$
 P_{\mu,\e} := -\frac{1}{2\pi\im}\oint_\Gamma ({\cal L}_{\mu,\e}-\lambda)^{-1} \di \lambda 
$$
is a well defined projection that commutes with $ {\cal L}_{\mu,\e} $, and  
is analytic in $ (\mu, \e)$.  The range $ \cV_{\mu,\e} := Rg( P_{\mu,\e}) $  
is a four dimensional symplectic subspace and 
the four eigenvalues of  $ \cL_{\mu,\e} $ inside 
$ \Gamma $ are the eigenvalues  of $ \cL_{\e,\mu} $ restricted to $ \cV_{\mu, \e} $. Additionally, 
since  $\cL_{\mu,\e}$ is a Hamiltonian and reversible operator, 
the projector $P_{\mu,\e}$ is  complex skew-Hamiltonian, 
namely $ (J P_{\mu,\e})^* = - JP_{\mu,\e} $, 
and reversibility preserving, meaning $ \bro P_{\mu,\e} = P_{\mu,\e}  \bro $.

These properties enable the construction of 
a ``symplectic and reversible basis'', analytic in $ (\mu,\e ) $,  
in which 
${\cal L}_{\mu,\e}\vert_{\cV_{\mu,\e}}$ is represented by a  $ 4\times 4$ 
Hamiltonian and reversible matrix 
\begin{equation}
\label{L.form}
\tL_{\mu,\e}= 
\begin{pmatrix} 
\tJ_2 E & \tJ_2 F \\ 
\tJ_2 F^* &\tJ_2 G 
\end{pmatrix} \, , 
\quad 
{\mathtt J}_2 := \begin{pmatrix} 
 0 & 1 \\
-1  & 0
\end{pmatrix} \, , \quad E = E^* \, , \ G = G^* \, ,  
\end{equation}
(the  matrices  $ E, G, F $ are  $ 2 \times 2 $) 
where  $\tJ_2 E$ has the 
expansion 
$$
\tJ_2 E =
- \im \big( \frac{\mu}{2}+ r(\mu\e^2,\mu^2\e,\mu^3) \big) +    
 \begin{pmatrix} 
0 & -\frac{\mu^2}{8}(1+r_5(\e,\mu))\\
- \e^2(1+r_1'(\e,\mu\e^2))+\frac{\mu^2}{8}(1+r_1''(\e,\mu))  & 
0 
 \end{pmatrix} \, . 
$$
The $ 2 \times 2 $ matrix $\tJ_2 E$  
exhibits the Benjamin-Feir phenomenon! It  
possesses two eigenvalues with non-zero real part as in 
\cref{eigelemu}, for $ 0 < \mu < \bar \mu (\e)  $ with $ \bar \mu (\e) \sim 2 \sqrt{2} \e $.  
In contrast,  the  $2\times 2$ 
matrix  $\tJ_2 G$ 
 has  purely imaginary eigenvalues for all $ \mu > 0 $, of order $\cO(\sqrt{\mu})$.

To prove that the complete $ 4 \times 4 $ matrix  
$ \tL_{\mu,\e} $  in \cref{L.form} exhibits the Benjamin-Feir phenomenon,   
we 
use a block diagonalization procedure inspired by KAM theory to eliminate the coupling term
 $  \tJ_2 F $. 
This  is a singular perturbation problem, since the spectra of the matrices 
$\tJ_2 E $ and $\tJ_2 G$ tend to $ 0 $ as $ \mu \to 0 $. 
However, by a Lie series expansion we can  conjugate the
matrix $\tL_{\mu,\e}$  in \cref{L.form} into a block diagonal  
 matrix 
$
  \begin{pmatrix} \mathtt{U} & \vline & 0 \\ \hline 0 & \vline & \mathtt{S} \end{pmatrix} $
 where  $ {\mathtt U } $ has the same form as $ \tJ_2 E  $
 and the  $ \mathtt{S} $ has purely imaginary 
 eigenvalues, as $ \tJ_2 G  $. This implies \cref{TeoremoneFinale}. 

\smallskip

In the finite-depth case in \cite{BMV3}, the four 
eigenvalues of $ \cL_{\mu,\e} $ near zero 
have all the  same size $ \cO( \mu ) $ (depending on $ \tth $) and more 
 steps in the Lie series expansion are required 
 to achieve an analoguous block diagonalization result. 
 Finally, for determining  
 the instability of the Stokes waves near the 
 Whitam-Benjamin critical depth 
$\tth_{{\mathtt W}{\mathtt B}}  $, 
  higher order expansions in $ (\e, \mu ) $ are
 needed, as shown in \cite{BMV_ed}. 
 
 \subsection{Unstable spectral bands away from zero.} 
The next  theorem establishes the existence of  infinitely many 
unstable spectral bands in the spectrum of $ \cL_{\e}$ away from zero, 
located close to  the  double eigenvalues 
$ \im \omega_*^{(\tp)}(\tth )  $ of $ \cL_{\umu,0}$ in \cref{eq:omegastar} 
 for sufficiently small values of 
 $ \e $ and Floquet exponents 
$ \mu $ near $ 	\umu $ defined in \cref{eq:umu}.

\begin{theorem}[\bf High frequency instabilities \cite{BCMV}] \label{Thm:1.3}
For any  integer   $ {\mathtt p} \geq 2$,   
there exist $\e_0>0$ 
and  a closed set of isolated depths  ${\cal S}^{({\mathtt p} )} \subset (0, + \infty)$ 
such that, for any $0 < \e \leq \e_0$, any $\tth \not\in 
{\cal S}^{({\mathtt p} )}$, the spectrum $\sigma_{L^2(\R)}(\cL_\e)$  contains 
${\mathtt p}-1$ disjoint ``isolas'' in the complex upper half-plane, plus the symmetric ones in 
the lower half-plane, see  \cref{fig:spectrum}, left. More precisely, 
there exists a non-zero analytic function  
\be\label{eq:zeri}
\beta_1^{(\tp)} :(0,+\infty) \to \R \, , \quad 
\ \tth \mapsto \beta_1^{(\tp)}(\tth) \, , \quad  \text{satisfying} \quad 
\lim_{\tth\to 0^+}
\beta_1^{(\tp)}
(\tth) = - \infty \, ,  \  \lim_{\tth\to+ \infty} \beta_1^{(\tp)} (\tth) = 0 \, , 
\ee
such that, for any $ \tth > 0 $
for which  
$ \beta_1^{(\tp)}(\tth) \neq 0 $,   the $\tp$-th unstable   isola  is approximated by an ellipse
\begin{enumerate}
\item[(i)] with 
semiaxis  having size $ \propto | \beta_1^{(\tp)}(\tth)| \e^\tp $;
\item[(ii)] 
centered  on the imaginary axis at a point  
$ \im y_0^{(\tp)}(\e) $ which is $ \cO(\e^2) $-close to $ \ii \omega_*^{(\tp)} (\tth) $,  cfr.  \cref{eq:omegastar}; 
\item[(iii)] 
 of  equation   $ x^2  + 
(E^{(\tp)}(\tth))^2(1+r(\e^2)) \big( y-y_0^{(\tp)}(\e) \big)^2 =   (\beta_1^{(\tp)}(\tth))^2 \e^{2\tp} (1+r(\e^{{2}})) $
where $ E^{(\tp)}(\tth)  \in (0,1) $.
\end{enumerate} 
 The  
  $\tp$-th unstable  isola 
 is traced by  a pair of eigenvalues $ \lambda_\tp^\pm (\mu,\e) $
  of $ \cL_{\mu,\e}$   
 which have  non zero real part 
iff  the  Floquet parameter 
$ \mu $  varies in a very narrow 
interval near $ \umu $ having width  $ \propto |\beta_1^{(\tp)}(\tth)|  \e^{\tp} $. 
\end{theorem}

Let us make some comments.  
\\[1mm]
1. {\sc  
Non degeneracy of the function $\beta_1^{(\tp)}(\tth) $.}
A challenging task  
 is to  prove the limits in \cref{eq:zeri} for {\it any} $\tp\geq 2 $, 
implying, in particular, that the  analytic function  
$\tth \mapsto \beta^{(\tp)}_1(\tth)$  is not identically  zero. The limit 
$|\beta_1^{(\tp)}(\tth)| \to + \infty   $ as $ \tth \to 0^+ 
$ shows that the high frequency instabilities become predominant in shallow water.  
\\[1mm]
2. {\sc Degenerate depths.
} 
If  $ \beta_1^{(\tp)} (\tth) = 0 $ 
the ellipses in \cref{Thm:1.3} shrink to a point, and 
the real part of the unstable eigenvalues of $ \cL_{\mu,\e} $ (if any) 
 is of much smaller size $\cO(\e^{\tp+2}) $. 
This degeneracy occurs at infinite depth  $ \tth = + \infty $ for {\em every} $\tp \geq 2 $,  
representing another aspect of  the integrability of the deep water waves equations.
 For $\tp = 2 $,   the first unstable isola  has been  proved  to exist at infinite depth in \cite{BMV4}, with  
 a semiaxis of size $ \propto \e^4  $.  

\subsection{Open problems.} We list a series of open problems 
 regarding the stability/instability of Stokes waves.
\begin{enumerate}
 \item[$P_{3,1}$)] {\it Large amplitude Stokes waves.}
What is the spectrum of the linearized operator $\cL_\e$ in \cref{cLepsilon}, both in finite and infinite depth,  as the amplitude  $\e$ of the Stokes  wave increases? Recent numerical simulations 
in \cite{DDLS} 
reveal a pletora of unstable spectral bands forming figures ``8'' and ``$ \infty $''. 

\item[$P_{3,2}$)] {\it Quasi-periodic  Stokes waves.}
Are the quasi-periodic traveling Stokes waves described in \cref{thm:main0}, and  \cite{BFM1,FG}, 
modulationally stable or unstable?  
Can we construct 
quasi-periodic unstable Bloch waves? 
\item[$P_{3,3}$)] {\it Instability of Stokes waves in 3D.}
Recently  
 \cite{CNS} proved that Stokes waves are linearly unstable under purely transversal perturbations,  and
  \cite{JRSY} 
under some longitudinal and transversal perturbations in finite depth.
The full shape of the spectrum under 3D spatial disturbances 
remains still unknown.

    \item[$P_{3,4}$)] {\it Nonlinear  Instability.} 
 Recently 
 \cite{ChenSu} proved that  small-amplitude Stokes waves in deep water are 
 nonlinearly unstable.
Is it possible to 
construct stable and unstable manifolds at the Stokes waves?
\end{enumerate}

\bibliographystyle{siamplain}



\end{document}